\documentclass[12pt]{article}
\usepackage{amsmath,amssymb,showlabels}

\newtheorem{Theorem}{Theorem}[part]
\newtheorem{Definition}{Definition}[part]

\newtheorem{Assumption}{Assumption}[part]

\newtheorem{Remark}{Remark}[part]

\def\cal#1{\mathcal{#1}}

\def\p{\partial}

\def\g{\gamma}

\def \F{I\!\!F}
\def \H{I\!\!H}

\def \N{I\!\!N}
\def \R{I\!\!R}

\def\Fc{{\cal F}}

\def\k{\kappa}

\def\T{{\cal T}}

\def\Dzw1#1{\frac{\partial^2 #1}{\partial z \partial w_1}}

\def\Dzb1#1{\frac{\partial^2 #1}{\partial z \partial b_1}}

\newcommand{\dproof}{\noindent {Proof.} \quad}
\newcommand{\fproof}{\hfill $\square$ \bigskip}

\newtheorem{definition}{Definition}[section]

\newtheorem{theorem}[definition]{Theorem}

\newtheorem{remark}[definition]{ \it Remark}

\newtheorem{proposition}[definition]{Proposition}
\newtheorem{lemma}[definition]{Lemma}

\def\BC{\mathcal{B}}

\def\PC{\mathcal{P}}

\def\R{{\bf R}}

\def\1B{\text{1\!\!I}}

\def\PC{\mathcal{P}}

\def\R{{\bf R}}

\def\1B{\text{1\!\!I}}

\def\p{\partial}

\begin{document}

\title{Generalized Dynkin Games and Doubly Reflected BSDEs with Jumps}
\author{Roxana Dumitrescu\thanks{CEREMADE,
Universit\'e Paris 9 Dauphine, CREST  and  INRIA Paris-Rocquencourt, email: {\tt roxana@ceremade.dauphine.fr}. The research leading to these results has received funding from the R\'egion Ile-de-France. }
\and
Marie-Claire Quenez\thanks{LPMA, 
Universit\'e Paris 7 Denis Diderot, Boite courrier 7012, 75251 Paris cedex 05, France, 
email: {\tt quenez@math.univ-paris-diderot.fr}}
\and 
Agn\`es Sulem
\thanks{INRIA Paris-Rocquencourt, Domaine de Voluceau, Rocquencourt, BP 105, Le Chesnay Cedex, 78153, France,  and Universit\'e Paris-Est, email: {\tt agnes.sulem@inria.fr}}}

\maketitle

\begin{abstract}
We introduce a generalized Dynkin game problem with non linear conditional expectation ${\cal E}$ induced by a Backward Stochastic Differential Equation (BSDE) with jumps. Let $\xi, \zeta$  be two RCLL adapted processes with 
$\xi \leq \zeta$. The criterium is given by 
\begin{equation*} 
 {\cal J}_{\tau, \sigma}= {\cal E}_{0, \tau \wedge \sigma } \left(\xi_{\tau}\textbf{1}_{\{ \tau \leq \sigma\}}+\zeta_{\sigma}\textbf{1}_{\{\sigma<\tau\}}\right)
 \end{equation*}
where $\tau$ and $ \sigma$ are stopping times valued in $[0,T]$. Under Mokobodski's condition, we establish the existence of a value function for this game, i.e. $\inf_{\sigma}\sup_{\tau} {\cal J}_{\tau, \sigma} = \sup_{\tau} \inf_{\sigma} {\cal J}_{\tau, \sigma}$. This value can be characterized via a doubly reflected BSDE. Using this characterization, we  provide some new results on these equations, such as comparison theorems and a priori estimates.
When $\xi$ and $\zeta$ are left upper semicontinuous along stopping times, we prove the existence of a saddle point. We also study  a generalized mixed game problem when the players have two actions: continuous control and stopping.
We then address the generalized Dynkin game in a Markovian framework and its links with parabolic partial integro-differential variational inequalities with two obstacles.
\end{abstract}

\vspace{10mm}

\noindent{\bf Key words~:}  Dynkin game, mixed Dynkin game,  $g$-expectation, non linear expectation, Backward stochastic differential equations with jumps, doubly reflected BSDEs, comparison theorem, partial integro-differential variational inequalities, viscosity solution.

\vspace{10mm} 

\noindent{\bf AMS 1991 subject classifications~:} 93E20, 60J60, 47N10.

\section{Introduction}

The classical Dynkin game  has been widely studied: see e.g. Bismut \cite{Bismuth}, 
Alario-Nazaret et al. \cite{ALM},  Kobylanski et al. \cite{KQC}.
Let $\xi, \zeta$  be  two  Right Continuous Left-Limited (RCLL) adapted processes with 
$\xi \leq \zeta$ and $\xi_T = \zeta_T$ a.s.\, The criterium is given, for each pair $(\tau, \sigma)$ of stopping times valued in $[0,T]$,  by 
\begin{equation*} 
 { J}_{\tau, \sigma}= {E} \left(\xi_{\tau}\textbf{1}_{\{ \tau \leq \sigma\}}+\zeta_{\sigma}\textbf{1}_{\{\sigma<\tau\}}\right).
 \end{equation*}

Under Mokobodski's condition,  which states that there exists two supermartingales such 
that their difference is between $\xi$ and  $\zeta$,  there exists a value function 
for the  Dynkin game, i.e. 
$\inf_{\sigma}\sup_{\tau} {J}_{\tau, \sigma} = \sup_{\tau} \inf_{\sigma} { J}_{\tau, \sigma}$.
When $\xi_t < \zeta_t$, $t <T$, and when $\xi$ and  $\zeta$ are also left upper semicontinuous, 
it is proved  that there exists a saddle point. 

Using a change of variable, these results can be generalized to the case of a criterium 
with  an instantaneous reward process $(g_t)$, of the form 
\begin{equation}\label{eq1}
 E  \left(\int_0^{\tau \wedge \sigma} g_s ds + \xi_{\tau}\textbf{1}_{\{ \tau \leq \sigma\}}+\zeta_{\sigma}\textbf{1}_{\{\sigma<\tau\}}\right).
\end{equation}
In the Brownian case and  when $(\xi_t)$ and $(\zeta_t)$ are 
continuous processes, Cvitani\'c  and Karatzas have established  in \cite{CK} 
links between these Dynkin games and doubly reflected  Backward stochastic differential equations with driver process $(g_t)$ and barriers $(\xi_t)$ and $(\zeta_t)$.

In this paper,  we introduce a generalization of the classical Dynkin game problem to the case of $g$-conditional expectations. Nonlinear expectations induced by BSDEs 
 have been introduced by S. Peng \cite{Peng} in the Brownian framework . Given a Lipschitz driver $g(t,y,z)$, a stopping time $\tau \leq T$ and a  square integrable ${\cal F}_S$-measurable 
 random variable $\eta$, the associated conditional $g$-expectation process denoted by $({\cal E}_{t,\tau}, 0 \leq t \leq \tau)$ is defined as the solution of the BSDE with driver $g$ and terminal 
 conditions $(\tau,\eta)$. The extension to the  case with jumps is  studied in \cite{R} and \cite{16}.
We consider the following generalized Dynkin game problem where the criterium is given, for each pair $(\tau, \sigma)$ of stopping times valued in $[0,T]$,  by 
\begin{equation*} 
 {\cal J}_{\tau, \sigma}= {\cal E}_{0, \tau \wedge \sigma } \left(\xi_{\tau}\textbf{1}_{\{ \tau \leq \sigma\}}+\zeta_{\sigma}\textbf{1}_{\{\sigma<\tau\}}\right)
 \end{equation*}
where  $\xi, \zeta$  are  two RCLL adapted processes with 
$\xi \leq \zeta$.

When  the driver $g$ does not depend on the solution, that is, when it is given by a process $(g_t)$, the criterium $ {\cal J}_{\tau, \sigma}$ coincides with \eqref{eq1}.
 It is well-known that  in this case, under Mokobodski's condition, the value function for the Dynkin game problem can be characterized as the solution of the Doubly Reflected BSDE (DRBSDE) associated  with driver process $(g_t)$ and barriers $(\xi_t)$ and $(\zeta_t)$ (see e.g. \cite{CK,HL,Lepeltier-Xu}). We generalize this result to the case of a non linear driver $g$ depending on the solution. More precisely, under Mokobodski's condition, we prove that  
  $$\inf_{\sigma}\sup_{\tau} { \cal J}_{\tau, \sigma} = \sup_{\tau} \inf_{\sigma} {\cal J}_{\tau, \sigma}$$
  and we characterize this common value function as  the solution of the DRBSDE associated  with driver $g$ and barriers $(\xi_t)$ and $(\zeta_t)$. Moreover, when $\xi$ and  $\zeta$ are left-upper semicontinuous along stopping times, we show that there exist saddle points. Note that, contrary to the previous existence results given in the case of classical Dynkin games, we do not assume the strict separability of the barriers.
  
   Then,  using the characterization of the solution of a DRBSDE as the value function of a generalized Dynkin game, we prove some results on DRBSDEs, such as a comparison 
   (respectively strict comparison) theorem and a priori estimates, which complete those given in the previous literature.
  

Moreover, we study  a generalized mixed game problem when the players have two actions: continuous control and stopping. The first (resp. second) player chooses a pair $(u, \tau)$ (resp. $(v, \sigma)$) of control and stopping time, and aims to maximize (resp. minimize) the criterium. In the previous literature (see \cite{BL} and \cite{HL}), the  criterium is given, for each quadruple $(u,\tau,v,\sigma)$ of controls and stopping times, by
\begin{equation}\label{uv}
  E_{Q^{u,v}}\left[\int_0^{\tau \wedge \sigma}c(t,u_t,v_t)dt+\xi_{\tau}\textbf{1}_{\{ \tau \leq \sigma\}}+\zeta_{\sigma}\textbf{1}_{\{\sigma<\tau\}}\right],
  \end{equation}
   where $Q^{u,v}$ are a priori  probability measures and $c(t,u_t,v_t)$ represents the instantenous reward associated with controls $u,v$. In this paper, we consider the following generalized mixed game problem. We are given a family of Lipchitz drivers $g^{u,v}$ and  the criterium is defined by
   \begin{equation}\label{nouveau}
\mathcal{E}_{0, \tau \wedge \sigma}^{u,v} \left(\xi_{\tau}\textbf{1}_{\{ \tau \leq \sigma\}}+\zeta_{\sigma}\textbf{1}_{\{\sigma<\tau\}}\right),
\end{equation} where $\mathcal{E}^{u,v}$ corresponds to the $g^{u,v}$-conditional expectation. Note that in the case of linear drivers
$g^{u,v}$\,, the criterium \eqref{nouveau} corresponds to 
a criterium of the form \eqref{uv}. In this particular case, when $\xi$ and $\zeta$ are regular, Hamad\`ene and Lepeltier have established some 
links between this mixed game problem and DRBSDEs (see \cite{HL}).
In this paper, we generalize these results to the case of non linear expectations and irregular payoffs $\xi$ and $\zeta$. We provide some sufficient  conditions which ensure the existence of a value function for the above generalized mixed game problem, and show that the common value function can be characterized as the solution of a DRBSDE. Under additional regularity assumptions on $\xi$ and $\zeta$, we show the existence of saddle points. 

 Finally, we address  the generalized Dynkin game in the Markovian case and its links with parabolic partial integro-differential variational inequalities (PIDVI) with two obstacles.

 The paper is organized as follows. In Section \ref{sec2} we introduce  notation and definitions and provide some preliminary results. In Section  \ref{sec3}, we consider a classical Dynkin game problem and study its links with
a DRBSDE associated with a driver which does not depend on the solution. We also provide an existence result for this game problem under relatively weak assumptions on $\xi$ and $\zeta$. In Section  \ref{sec4}, we introduce a generalized Dynkin game problem expressed in terms of  $g$-conditional expectations. We prove the existence of a value function for this game problem. We show that the common value function can be characterized as the solution of a non linear DRBSDE with jumps and RCLL barriers $\xi$ and $\zeta$.  We then study  a generalized mixed game problem when the players have two actions: continuous control and stopping.  In Section  \ref{sec5}, we provide comparison theorems   and a priori estimates for DRBSDEs with jumps and RCLL obstacles.  In  the Markovian case,  relations between generalized Dynkin games  and PIDVIs are studied in Section  \ref{sec6}. We state that the value function of the generalized Dynkin game corresponds to a solution of a PIDVI in the viscosity sense. Under additional assumptions, we obtain an uniqueness result in the class of continuous and bounded functions.

\section{Notation and definitions} \label{sec2}

Let $(\Omega,  \Fc_, P)$ be a probability space. 
Let  $W$ be   a  one-dimensional Brownian motion.  Let $({\bf E}, {\cal K})$  be a measurable space 
equipped with a $\sigma$-finite positive measure $\nu$ and
 let  $N(dt,de)$ be a Poisson random measure with compensator $\nu(de)dt$.
 Let $\tilde N(dt,de)$ be its compensated process.
Let  $\F = \{\Fc_t , t \geq 0 \}$ 
be  the natural filtration associated with $W$ and $N$. 

\paragraph{Notation.}

Let ${\cal P}$ be  the predictable $\sigma$-algebra
on $[0,T]  \times \Omega$.

 For each $T>0$, we use the following notation:  $L^2({\cal F}_T)$  is the set of random variables $\xi$ which are  $\Fc
_T$-measurable and square integrable; $\H^{2}$ is the set of
real-valued predictable processes $\phi$ such that $\| \phi\|^2_{\H^{2}} := E \left[\int_0 ^T \phi_t ^2 dt \right] < \infty$; ${\cal S}^{2}$ denotes the set of real-valued RCLL adapted
 processes $\phi$ such that
$\| \phi\|^2_{{\cal S}^2} := E(\sup_{0\leq t \leq T} |\phi_t |^2) <  \infty$; ${\cal A}^2$ (resp. ${\cal A}^1$) is the set of real-valued non decreasing RCLL predictable
 processes $A$ with $A_0 = 0$ and $E(A^2_T) < \infty$ (resp. $E(A_T) < \infty$). We also introduce the following spaces.

\begin{itemize}
%
%
\item $L^2_\nu$ is the set of Borelian functions $\ell: {\bf E} \rightarrow \R$ such that  $\int_{{\bf E}}  |\ell(e) |^2 \nu(de) < + \infty.$

The set  $L^2_\nu$ is a
 Hilbert space equipped with the scalar product
$$\langle \delta    , \, \ell \rangle_\nu := \int_{{\bf E}} \delta(e) \ell(e) \nu(de)  \quad \text{ for all  }  \delta  , \, \ell \in L^2_\nu \times L^2_\nu,$$
and the norm $\|\ell\|_\nu^2 :=\int_{{\bf E}}  |\ell(e) |^2 \nu(de).$


\item $\H_{\nu}^{2}$ is  the set of processes $l$ which are {\em predictable}, that is, measurable $$l : ([0,T]  \times \Omega \times  {\bf E},\; \PC \otimes {\cal K}  \rightarrow (\R\;,  \BC(\R)); \quad
(\omega,t,e) \mapsto l_t(\omega, e)
$$
such that $\| l \|^2_{\H_{\nu}^{2}} :=E\left[\int_0 ^T \|l_t\|_{\nu}^2 \,dt \right]< \infty.$





\end {itemize}

Moreover,  $\T_{0}$ is the set of
stopping times $\tau$ such that $\tau \in [0,T]$ a.s.\, and for each $S$ in $\T_{0}$, we denote by 
   $\T_{S}$ the set of
stopping times
$\tau$ such that $S \leq \tau \leq T$ a.s.



\begin{definition}[Driver, Lipschitz driver]\label{defd}
A function $g$ is said to be a {\em driver} if
\begin{itemize}
\item
$g: [0,T]  \times \Omega \times \R^2 \times L^2_\nu \rightarrow \R $\\
$(\omega, t,y, z, \k(\cdot)) \mapsto  g(\omega, t,y,z,k(\cdot)) $
  is $ {\cal P} \otimes {\cal B}(\R^2)  \otimes {\cal B}(L^2_\nu)
- $ measurable,
\item $g(.,0,0,0) \in \H^2$.
\end{itemize}
A driver $g$ is called a {\em Lipschitz driver} if moreover there exists a constant $ C \geq 0$ such that $dP \otimes dt$-a.s.\,,
for each $(y_1, z_1, k_1)$, $(y_2, z_2, k_2)$,
$$|g(\omega, t, y_1, z_1, k_1) - g(\omega, t, y_2, z_2, k_2)| \leq
C (|y_1 - y_2| + |z_1 - z_2| +   \|k_1 - k_2 \|_\nu).$$
\end{definition}


%

%
Recall that for each Lipschitz driver $g$, and each terminal condition $\xi$ $\in$ $L^2({\cal F}_T)$, there exists a unique solution $(X, \pi, l)$ $\in$ ${\cal S}^{2} \times \H^{2} \times \H_{\nu}^{2}$ satisfying
 \begin{equation}\label{er}
-dX_t =g(t,X_{t^-},\pi_t, l_t(\cdot))dt - \pi_t dW_t - \int_{{\bf E}}  l_t(e) \tilde N(dt,de)  ; \qquad  X_T = \xi.
\end{equation}

 The  solution is denoted by $(X( \xi,T), \pi ( \xi,T), l( \xi,T))$.

This result can be extended when the terminal time  $T$ is replaced by a
stopping time $\tau \in \T_0$ and when $\xi$ is replaced by a random variable $\eta \in L^2({\cal F}_\tau)$. The solution $X_\cdot(\eta, \tau)$ corresponds to the so-called $g$-conditional expectation of $\eta$, denoted by $\mathcal{E}_{\cdot,\tau}(\eta).$

\begin{definition}\label{proba}
 Let $A=(A_t)_{0 \leq t \leq T}$ and $A'=(A_t')_{0 \leq t \leq T}$ belonging to $\cal{A}^1$. We say that the measures $dA_t$ and $dA'_t$ are {\em mutually singular}, and we write $dA_t \perp dA'_t$\,, if there exist $D \in \mathcal{P}$ such that:
$$\int_0^T \textbf{1}_{D^c_t} dA_t= \int_0^T \textbf{1}_{D_t} dA'_t=0 \,\, \text{ a.s.}\,,  $$
where for each $t \in [0,T]$, $D_t$ is the {\em section at time $t$} of $D$, that is, $D_t:= \{\omega \in \Omega \,,\, (t, \omega) \in D \}$.
\end{definition}

We introduce the DRBSDEs with jumps, for which the solution is constrained to stay between two given RCLL processes called barriers $\xi \leq \zeta$. Two nondecreasing processes A and $A'$ are introduced in order to push the solution $Y$ above $\xi$ and below $\zeta$, and this in a minimal way. This \textit{minimality property} of $A$ and $A'$ is ensured by  the so called Skorohod conditions (see condition $(iii)$ below) together with the additional constraint $dA_t \perp dA'_t$\\     (see condition $(ii)$ below). 
\begin{definition}[Doubly Reflected BSDEs with Jumps]
Let $T>0$ be a fixed terminal time and  $g$ be  a Lipschitz driver.doubly 
Let $\xi$  and $\zeta$ be two adapted RCLL processes with $\zeta_T= \xi_T$ a.s.,   $\xi \in {\cal S}^2$,  $\zeta \in {\cal S}^2$, $\xi_t \leq \zeta_t$, $0 \leq t \leq T$ a.s.

A process $(Y, Z, k(.), A,A')$  in ${\cal S}^2 \times \H^2 \times \H^2_{\nu}\times {\cal A}^2 \times {\cal A}^2$  is said to be a solution of the doubly  reflected BSDE (DRBSDE)  associated with driver $g$ and barriers $\xi, \zeta$ if
 \begin{align}
   -dY_t &= g(t,Y_t,  Z_t, k_t(\cdot) )dt +dA_t -dA_t^{'}- Z_t  dW_t -\int_{{\bf E}}  k_t(e) \tilde{N}(dt,de); \;  Y_T = \xi_T, \label{DRBSDE} \\
   \text{with} & \nonumber \\
  &(i)   \;\; \xi_t \leq Y_t \leq \zeta_t,  \; 0 \leq t \leq T \text{ a.s.}, \nonumber \\
 &(ii)   \; \; dA_t \perp dA_t'  \nonumber \\
& (iii)   \displaystyle   \int_0^T (Y_t - \xi_t) dA^c_t = 0 \text{ a.s. and } \; \displaystyle    \int_0^T (\zeta_t-Y_t)dA^{'c}_t = 0 \text{ a.s.  } \nonumber \\
& \qquad  \Delta A_{\tau}^d=  \Delta A_{\tau}^d {\bf 1}_{\{Y_{\tau^-} = \xi_{\tau^-}\}} \text{ and }  \;  \Delta A_{\tau}^{'d}= \Delta A_{\tau}^{'d} {\bf 1}_{\{Y_{\tau^-} = \zeta_{\tau^-}\}} \text{ a.s. } \forall \tau \in \T_0 \text{ predictable } \nonumber
\end{align}
\end{definition}

Here  $A^c$ (resp $A^{'c}$) denotes the continuous part of $A$ (resp $A^{'}$) and $A^d$ (resp $A^{'d}$) its discontinuous part.

\begin{remark}
The above definition is not exactly the same as the one given in the previous literature, where $A$ and $A'$ are not constrained to satisfy $dA_t \perp dA_t'$. Note that when $A$ and $A'$ are not required to be mutually singular, they can simultaneously increase 
on $\{\xi_{t^-}=\zeta_{t^-}\}$.
\end{remark}

We  introduce the following definition.

\begin{Definition} \label{defr} A progressive process $(\phi_t)$ (resp. integrable) is said to be {\em left-upper semicontinuous (l.u.s.c.) along stopping times} (resp. along stopping times in expectation ) if for all $\tau \in \T_0$ and for each non decreasing sequence of stopping times $ (\tau_n)$ such that $\tau^n \uparrow \tau$ a.s.\,,
\begin{equation}\label{usc}
\phi_{\tau} \geq \limsup_{n\to \infty} \phi_{\tau_n} \quad \mbox{a.s.} \quad \text{ (resp. } E[\phi_{\tau}] \geq \limsup_{n\to \infty} E[\phi_{\tau_n}] ).
\end{equation}
\end{Definition}

\begin{remark}\label{lusc}
Note that  when $(\phi_t)$ is left-limited, then $(\phi_t)$ is {\em left-upper semicontinuous (l.u.s.c.) along stopping times} if and only if for all predictable stopping time $\tau \in \T_0$, 
$\phi_{\tau} \geq  \phi_{\tau_-}\quad\mbox{a.s.}$
\end{remark}


\section{Classical Dynkin games and links with doubly reflected BSDEs with a driver process  }\label{sec3}

In this section, we are given a predictable process $g= (g_t)$ in $\mathbb{H}^2$.\\
Let $\xi$  and $\zeta$ be two adapted  processes only supposed to be RCLL with $\zeta_T= \xi_T$ a.s.,   $\xi \in {\cal S}^2$,  $\zeta \in {\cal S}^2$, $\xi_t \leq \zeta_t$, $0 \leq t \leq T$ a.s. \\
We state that 
the doubly reflected BSDE associated with the driver process $(g_t)$ and the barriers  $\xi$  and $\zeta$ admits a  unique solution $(Y,Z,k(\cdot),A,A')$, which is related to a classical Dynkin game problem defined below. Our results complete the previous works on classical Dynkin games and DRBSDEs ( see for e.g. \cite{CK}, \cite{HH}). In particular, we provide an existence result of saddle points under weaker assumptions than those made in the previous literature.\\
For any $S \in \T_{0}$ and any stopping times $\tau, \sigma \in \T_S$, consider the gain (or payoff):
\begin{equation} \label{payoff}
 I_S(\tau, \sigma)=\int_{S}^{\sigma \wedge \tau}g(u)du+\xi_{\tau}\textbf{1}_{\{ \tau \leq \sigma\}}+\zeta_{\sigma}\textbf{1}_{\{\sigma<\tau\}}.
 \end{equation}
For any $S \in \T_{0}$, the upper and lower value functions at time $S$ are defined respectively by 
\begin{equation} \label{overline}
\overline{V}(S) :=ess \inf_{\sigma \in \T_S } ess \sup_{\tau \in \T_S} E[I_S( \tau, \sigma)| \cal{F}_{S}]
  \end{equation}

\begin{equation} \label{underline}
 \underline{V}(S):=ess\sup_{\tau\in \T_S } ess \inf_{\sigma \in \T_S} E[I_S(\tau, \sigma)| \cal{F}_{S}]
 \end{equation}
 We clearly have the inequality  $\underline{V}(S) \leq \overline{V}(S)$ a.s.\\
By definition, we say that there {\em exists a value function} at time $S$ for the Dynkin game problem if $\overline{V}(S)=\underline{V}(S)$ a.s.
\begin{definition}[$S$-saddle point]
Let $S \in \T_0$. A pair $(\tau^*, \sigma^*) \in \T_S^2$ is called an {\em $S$-saddle point} if for each $(\tau, \sigma) \in \T_S^2$, we have
$$  E[I_S(\tau, \sigma^*)| \cal{F}_{S}]
 \leq  E[I_S(\tau^*, \sigma^*)| \cal{F}_{S}] \leq E[I_S(\tau^*, \sigma)| \cal{F}_{S}] \;\; \text {a.s. }$$
\end{definition}
We introduce the following RCLL adapted processes which depend on the process $g$:
\begin{equation}\label{defz}
\Tilde{\xi}_t^{g}:=\xi_t-E[\xi_T+\int_t^Tg(s)ds|\cal{F}_t], \;
\quad \Tilde{\zeta}_t^{g}:=\zeta_t-E[\xi_T+\int_t^Tg(s)ds|\cal{F}_t], \quad 0 \leq t \leq T.
\end{equation}
They satisfy the important property 
$\Tilde{\xi}_T^{g} = \Tilde{\zeta}_T^{g} = 0 \; \text{ a.s. }$
Moreover, this change of variables allows us to get rid of the term $\int g(t)dt$, and thus to simplify the notation. Some more comments on this change of variables are given in Remark $\ref{8}$ in the Appendix.

For each RCLL adapted process  $\phi = (\phi_t)_{0 \leq t \leq T} $ valued in $\R \cup \{+\infty \}$ with $\phi^{-} \in \mathcal{S}^2$, we denote by $\mathcal{R}(\phi)$ the Snell envelope of $\phi$, defined as the minimal RCLL supermartingale greater or equal to $\phi$ a.s.\, By the optimal stopping theory, $\mathcal{R}(\phi)$ is equal to the value function of the optimal stopping problem associated with the reward $\phi$.

We state the following lemma. 
\begin{lemma} \label{seq} 
There exists a  unique pair of non-negative RCLL supermartingales $(J^g, J'^g)$ valued in 
$[0, +\infty]$ satisfying $J^g_T=J'^g_T=0$ a.s. and the system
 \begin{equation}\label{sys2}
J^g={\cal R}(J'^g+\Tilde{\xi}^g) \quad ; \quad  J'^g={\cal R}(J^g-\Tilde{\zeta}^g).
\end{equation}
and satisfying the following minimality property: if $H$ and $H'$ are non-negative RCLL supermartingales valued in 
$[0, +\infty]$
such that $H \geq H'+\tilde{\xi}^g $ and  $H' \geq H - \tilde{\zeta}^g$ , then we have $J^g\leq H$ and $J^{'g} \leq H'$.

\end{lemma}


A sketch of the proof is given in the Appendix. Using this lemma, we derive the following result.

\begin{theorem}\label{OY}
Let $\xi$  and $\zeta$ be two adapted RCLL processes in ${\cal S}^2$ with $\zeta_T= \xi_T$ a.s.\,and  $\xi_t \leq \zeta_t$, $0 \leq t \leq T$ a.s. 
Suppose that  $J^g, J^{\prime \, g}\in {\cal S}^2$.\\
Let $\overline Y$ be the RCLL adapted process defined by  
\begin{equation}\label{oY}
{\overline Y}_t:= J^g_t-J^{'g}_t+E[\xi_T+\int_t^Tg(s)ds | \cal{F}_t] ; \; 0 \leq t \leq T.
\end{equation}
 There exist $(Z, k, A, A' ) \in \H^2 \times \H^2_{\nu}\times {\cal A}^2 \times {\cal A}^2$  such that $({\overline Y}, Z, k, A,A')$  is  a solution of DRBSDE \eqref{DRBSDE} 
 associated with the driver process $g(t)$.
 \end{theorem}

 \dproof
By assumption, $J^g$ and $J^{'g}$  are square integrable supermartingales. The process $\overline Y$ is thus well defined. 
By Lemma \ref{seq}, we have $J^g_T=  J'^g_T$ a.s.   Hence, $\overline Y_T = \xi_T$ a.s.\,
By the Doob-Meyer decomposition, there exist two square integrable martingales $M$ and $M'$ and two processes $B$ and $B^{'}$ $\in$ ${\cal A}^2$ such that:
\begin{equation}\label{1}
dJ_t^g=dM_t-dB_t \quad ; \quad dJ_t^{'g}=dM_t^{'}-dB_t^{'}.
\end{equation}
 Set 
 $$\overline{M}_t :=M_t-M_t^{'}+E[\xi_T+\int_0^T g(s)ds | \cal{F}_t].$$ 
 By \eqref{1},  \eqref{oY},   we derive $d\overline Y_t=
d\overline{M}_t-d\alpha_t -g(t)dt$, with $\alpha:= B-B^{'}$. Now, by the martingale representation theorem, there exist $Z \in \mathbb{H}^{2}, k \in \mathbb{H}_{\nu}^{2}$ such that
$
d\overline{M}_t=Z_tdW_t+\int_{{\bf E}} k_t(e)\Tilde{N}(de,dt)
$.
Hence, $$ -d\overline Y_t = g(t )dt +d\alpha_t- Z_t  dW_t -\int_{{\bf E}}  k_t(e) \tilde{N}(dt,de).$$

Let us now show that $B,B'$ satisfy the Skorohod conditions \eqref{DRBSDE}(iii).\\
 By the optimal stopping theory (see e.g. Proposition B.1 in  \cite{3}), the process $B^c$ increases only when the value function $J^g$ is equal to the corresponding reward $J^{'g}+\Tilde \xi^g$. Now, $\{J_t^g=J_t^{'g}+\Tilde \xi^g\}=\{ \overline{Y_t}= \xi_t\}$. Hence, $\int_{0}^{T}(\overline{Y_t}-\xi_t)dB_t^c=0$ a.s.  Similarly the process $B^{'c}$ satisfies $\int_{0}^{T}(\overline{Y_t}-\zeta_t)dB_t^{'c}=0$ a.s.\\
  Moreover, for each predictable stopping time $\tau \in \T_0$ we have 
$\Delta{B_{\tau}^d}=\textbf{1}_{J_{\tau^-}^g=J_{\tau^-}^{'g}+\Tilde \xi_{\tau^-}^g}\Delta{B_{\tau}^d}=\textbf{1}_{\overline{Y}_{\tau^-}= \xi_{\tau^-}}\Delta B_\tau^d$ a.s. and $\Delta{B_{\tau}^{'d}}=\textbf{1}_{\overline{Y}_{\tau^-}= \zeta_{\tau^-}}\Delta B_\tau^{' d}$ a.s.

By Proposition \ref{canonique} in the Appendix, there exist $A, A' \in \mathcal{A}^2 $ such that $\alpha=A-A'$ with $dA_t \perp dA'_t$. Also, $dA_t << dB_t$. Hence, 
since $\int_0^T \textbf{1}_{Y_{t^-} > \xi_{t^-}}dB_t = 0$ a.s.\,, we get
 $ \int_0^T \textbf{1}_{Y_{t^-} > \xi_{t^-}}dA_t=0$ a.s.\, Similarly, we obtain $\int_0^T \textbf{1}_{Y_{t^-} < \zeta_{t^-}}dA'_t=0$ a.s. The processes $A$ and $A'$ thus satisfy the Skorohod conditions \eqref{DRBSDE}(iii).
\fproof

\begin{remark}
Set $H_t:=E[A_T-A_t | \mathcal{F}_t]$ (resp. $H'_t:=E[A'_T-A'_t | \mathcal{F}_t]$). Using the same notation as in the above proof, since $dA_t << dB_t$ (resp. $dA'_t << dB'_t$), we have $H_t \leq J_t= E[B_T-B_t | \mathcal{F}_t]$ (resp. $H'_t \leq J'_t
= E[B'_T-B'_t | \mathcal{F}_t]$). Moreover, $H-H'= J-J'$. Hence, we have $H \geq H'+\tilde{\xi}^g $ and  $H' \geq H - \tilde{\zeta}^g$. By the minimality property of $J$, $J'$ (see  Lemma \ref{seq}), we derive that $J_t=H_t=E[A_T-A_t | \mathcal{F}_t]$ (resp. $J'_t= H'_t= E[A'_T-A'_t | \mathcal{F}_t]$). 
\end{remark}

From this theorem, we derive the following uniqueness and existence result for the DRBSDE associated with the driver process $(g_t)$, as well as the characterization of the solution as the value function of the 
above Dynkin game problem. We also show that if the associated non decreasing processes $A$ and $A'$ are continuous, 
then there exist saddle points for this game problem.
\begin{theorem}\label{f}
Let $\xi$  and $\zeta$ be two adapted RCLL processes in ${\cal S}^2$ with $\zeta_T= \xi_T$ a.s.\,and $\xi_t \leq \zeta_t$, $0 \leq t \leq T$ a.s.\,Suppose that $J^g_t, J'^g_t\in \mathcal{S}^2$. 
 
The doubly reflected  BSDE~\eqref{DRBSDE} associated with driver process $g(t)$ admits a unique solution
 $(Y,Z,k,A,A')$ in ${\cal S}^2 \times \H^2 \times \H^2_{\nu}\times ({\cal A^2})^2$.\\
For each $S \in \T_0$, $Y_S$ is the common value function of the Dynkin game, that is
\begin{eqnarray}\label{deux-2}
 Y_S= \overline{V}(S)=\underline{V}(S) \quad \mbox {a.s.}
\end{eqnarray}
Moreover, if the processes $A,A'$ are continuous, then, for each $S$ $\in$ $\T_0$, the pair of stopping times $(\tau_s^*, \sigma_s^*)$ defined by 
\begin{equation}\label{qqq}
\sigma^*_S:= \inf \{ t \geq S,\,\, Y_t = \zeta_t\}; \quad \tau^*_S:= \inf \{ t \geq S,\,\, Y_t = \xi_t\}.
\end{equation}
is an $S$-saddle point for the Dynkin game problem associated with the gain $I_S$.

%
\end{theorem}

\begin{remark}\label{sauts}
For each predictable stopping time $\tau \in \T_0$, we have
$\Delta A_\tau^d=(\Delta Y_\tau)^- $ and $\Delta A_\tau^{'d}=  (\Delta Y_\tau)^+ \text{ a.s.}$
\end{remark}
A short proof is given in the Appendix. Note that the uniqueness of the non decreasing RCLL processes $A$ and $A'$ holds because of the constraint $dA_t \perp dA'_t$\, (see the Appendix for details).


We now provide a sufficient condition on $\xi$ and $\zeta$ for the existence of saddle points. By the last assertion of Theorem \ref{f}, 
it is sufficient to give a condition which ensures the continuity of $A$ and $A'$.
\begin{theorem}[Existence of $S$-saddle points]\label{continuous}
Suppose that the assumptions of Th. \ref{f} are satisfied and that  $\xi$ and $-\zeta$ are l.u.s.c. along stopping times.\\
 Let $(Y,Z,k(.),A,A')$ be the solution of DRBSDE \eqref{DRBSDE}.\\ 
The processes $A$ and $A'$ are then continuous. Also, for each $S \in \T_0$, the pair of stopping times $(\tau^*_S, \sigma^*_S )$ defined by \eqref{qqq}
 is an $S$-saddle point.
\end{theorem}

\begin{remark} The assumptions made on $\xi$ and $\zeta$ are milder than the ones made in the literature where it is also supposed $\xi_t<\zeta_t, t<T$ a.s. ( see  e.g. \cite{ALM}, \cite{CK}, \cite{KQC}).
\end{remark}

\dproof
By the second assertion of Theorem \ref{f}, it is sufficient to prove that $A$ and $A'$ are continuous.
Let $\tau \in \T_0$ be a predictable stopping time. Let us show $\Delta A_\tau=0$ a.s.\\
By Remark \ref{sauts}, we have $\Delta A_\tau=(\Delta Y_\tau)^- $ a.s.\\
Since $dA_t \perp dA'_t$ \,, there exists $D \in \mathcal{P}$ such that:
$\int_0^T \textbf{1}_{D_t^c} dA_t= \int_0^T \textbf{1}_{D_t} dA'_t=0 \,\, \text{ a.s.}  $ We introduce the set $D_\tau:= \{\omega,\, (\tau(\omega), \omega) \in D \}.$ Since $A$ satisfies the Skorohod condition, we have $D_\tau \subset \{Y_{\tau^-}=\xi_{\tau^-}\}$ a.s.
We thus obtain
$$
\Delta A_{\tau}=\bf{1}_{D_\tau}( Y_{\tau^-}-Y_{\tau})^+=\bf{1}_{D_\tau}( \xi_{\tau^-}-Y_{\tau})^+ \leq \bf{1}_{D_\tau}( \xi_{\tau}-Y_{\tau}) \quad {\rm a.s.}
$$
The last inequality follows from the inequality $\xi_{\tau^-} \leq \xi_{\tau}$ a.s. (see Remark \ref{lusc}). Since $ \xi \leq Y$, we derive that $\Delta A_{\tau} \leq 0$ a.s. Hence, $\Delta A_{\tau} = 0$ a.s.\,, and this holds for each predictable stopping time $\tau$. Consequently, $A$ is continuous. Similarly, one can show that $A'$ is continuous.
\fproof

 Since $J^g \geq J^{'g}+\Tilde{\xi}^g$ and $J^{'g} \geq J^{g}-\Tilde{\zeta}^g$, the condition $J^g \in {\cal S}^2$ is equivalent to the condition $J^{'g} \in {\cal S}^2$.

We now recall the definition of Mokobodski's condition.
\begin{definition}[Mokobodski's condition]
Let $\zeta, \xi \in \mathcal{S}^2$. The Mokobodski's condition is defined as follows:
there exist two nonnegative RCLL supermartingales $H$ and $H'$ $\in \mathcal{S}^2$ such that:
\begin{equation}\label{Moki}
\xi_t {\bf 1}_{t <T} \leq H_t -H'_t \leq \zeta_t {\bf 1}_{t <T} \quad 0\leq t \leq T \quad {\rm a.s.}.
\end{equation}
\end{definition}

\begin{proposition}\label{Mokoko}
Let $g \in \H^2$.  The following assertions are equivalent:
\begin{itemize}
\item[(i)] $J^g \in \mathcal{S}^2$ 
\item[(ii)]  $J^0 \in \mathcal{S}^2$
\item[(iii)]  Mokobodski's condition holds.
\item[(iv)]  DRBSDE \eqref{DRBSDE} with driver process $(g_t)$ has a  solution. 
\end{itemize}
\end{proposition}

\dproof
Using the minimality property of $J$ and $J'$ given in Lemma \ref{seq}, one can show that 
 $J^g \in {\cal S}^2$ if and only if  
  there exist two non-negative  supermartingales
$H^g,H^{'g} \in \mathcal{S}^2 $ such that 
\begin{equation}\label{Mokobis}
  \tilde \xi_t^g \leq H_t^g -H^{'g}_t \leq \tilde \zeta_t^g  \quad 0\leq t \leq T \quad {\rm a.s.}
\end{equation}
Since this equivalence  holds for all $g \in \H^2$, in particular  when $g =0$, we get (ii) $\Leftrightarrow$  (iii).
It remains to show  (i) $\Leftrightarrow$  (ii)
For this, it is sufficient to show that $\eqref{Moki}$ is equivalent to $\eqref{Mokobis}$.
Suppose that $\eqref{Moki}$ is satisfied. 
By setting
\begin{equation*}
\begin{cases}
H_t^{g}:=H_t-E[\xi_T^+(s)ds|\cal{F}_t]-E[\int_t^Tg^+(s)ds|\cal{F}_t], 0 \leq t \leq T\\
H_t^{'g}:=H'_t-E[\xi_T^-(s)ds|\cal{F}_t]-E[\int_t^Tg^-(s)ds|\cal{F}_t], 0 \leq t \leq T,\\
\end{cases}
\end{equation*}
$\eqref{Mokobis}$ holds. 
It remains to prove that $(iv)$ implies $(i)$. Let $(Y,Z,k,A,A')$ be the solution of the DRBSDE \eqref{DRBSDE} associated with driver process $(g_t)$.  Let $H_t^{g}:=E[A_T-A_t|\cal{F}_t]$ and $H_t^{'g}:=E[A_T^{'}-A_t^{'}|\cal{F}_t]$. We have $H_t^g-H_t^{'g}=Y_t-E[\int_t^Tg(s)ds|\cal{F}_t].$ Since $\xi \leq Y \leq \zeta$, condition \eqref{Mokobis} holds.
\fproof

\section{Generalized Dynkin games and links with doubly reflected BSDEs with a non linear driver}\label{sec4}
In this section, we are given a Lipschitz driver $g$.
\subsection{Existence and uniqueness for DRBSDEs}
\begin{theorem}\label{exiuni}
Suppose $\xi$ and $\zeta$ are RCLL adapted process in ${\cal S}^2$ such that $\xi_T = \zeta_T$ a.s.\,and $\xi_t \leq \zeta_t$, $0 \leq t \leq T$ a.s.\,Suppose that $J^0 \in {\cal S}^2$ (or equivalently suppose that Mokobodski's condition is satisfied).\\
Then, DRBSDE (\ref{DRBSDE}) admits a unique solution $(Y,Z,k(.), A,A')\in {\cal S}^2 \times \H^2 \times \H^2_{\nu}\times ({\cal A^2})^2$. \\    
If $\xi$ and $\zeta$ are  l.u.s.c. along stopping times, then the processes $A$ and $A'$ are continuous.
\end{theorem}
The proof is based on classical arguments. It is given in the appendix.

\begin{remark}

Note  that for each predictable stopping time $\tau \in \T_0$, we have
$\Delta A_\tau^d=(\Delta Y_\tau)^- $ and $\Delta A_\tau^{'d}=  (\Delta Y_\tau)^+ \text{ a.s.}$
\end{remark}

\subsection{Generalized Dynkin games}\label{caract}

In this section, we introduce a generalized Dynkin game expressed in terms of $g$-conditional expectations.\\
In order to ensure that the $g$-conditional expectation $\mathcal{E}$ is non decreasing, we make the following assumption.

\begin{Assumption}\label{Royer}
A  lipschitz driver $g$ is said to satisfy Assumption~\ref{Royer} if 
for each  process $(Y,Z,l^1,l^2)$ in $S^{2} \times \H^{2} \times (\H_{\nu}^{2})^2$, 
there exists a bounded predictable process $(\gamma_t)$ such that $dt\otimes dP \otimes \nu(du)$-a.s.\,,
\begin{equation}\label{robis}
\gamma_t(u) \geq -1  \;\; \text{ and }
\;\; |\gamma_t(u) | \leq \psi (u),
\end{equation}
where $\psi$ $\in$ $L^2_{\nu}$, and such that 
\begin{equation}\label{autre}
g(t, Y_t, Z_t, l^1_t) - g(t,Y_t,Z_t,l^2_t) \geq \langle \gamma_t\,,\, l^1_t- l^2_t \rangle_\nu , \;\; t \in [0,T],\; \; dt\otimes dP \text{ a.s.}.
\end{equation}

\end{Assumption}
 For example, this assumption is satisfied if $g$ is ${\cal C}^1$ with respect to $\ell$ with $\nabla_{\ell} g \geq -1$ and $\vert \nabla_{\ell} g  \vert \leq\psi $, where $\psi$ $\in$ $L^2_{\nu}$ (see Lemma \ref{diff} in the Appendix).\\
Moreover, the above assumption ensures the non decreasing property of $\mathcal{E}^g$ by the comparison theorem for BSDEs with jumps (see Theorem 4.2 in \cite{16}). In the case when in \eqref{robis}, we have $\gamma_t>-1$, by the strict comparison  theorem (see Theorem 4.4 in \cite{16}), it follows that $\mathcal{E}^g$ is strictly monotonous.

We now introduce the following game problem, which can be seen as a Dynkin game written in terms of $g$-conditional expectations.

For each $\tau, \sigma \in \T_0$, the {\em reward} at time $\tau \wedge \sigma$ is given by the random variable 
\begin{equation} \label{terminal}
I(\tau,\sigma):=\xi_{\tau} {\bf 1} _{\tau \leq \sigma}+ \zeta_{\sigma} {\bf 1} _{\sigma < \tau}.
\end{equation}
Note that $I(\tau, \sigma)$ is ${\cal F}_{\tau \wedge \sigma}$-measurable.

Let $S \in \T_0$. For each $\tau \in \T_S$  and $\sigma \in \T_S $, the associated {\em criterium} is given by $\cal{E}_{S,\tau \wedge \sigma}(I(\tau, \sigma))$, the $g$-conditional expectation of the reward $I(\tau, \sigma)$. 
Recall that $\cal{E}_{\cdot,\tau \wedge \sigma}(I(\tau, \sigma))= X_{\cdot}^{\tau, \sigma},$ where $(X_{\cdot}^{\tau, \sigma},
\pi_{\cdot}^{\tau, \sigma}, l_\cdot^{\tau, \sigma})$ is the solution of the BSDE associated with driver $g$, terminal time $\tau \wedge \sigma$ and terminal condition $I(\tau, \sigma)$, that is
\begin{eqnarray*}
 -dX_{s}^{\tau, \sigma}= g(s,X_{s}^{\tau, \sigma}, \pi_{s}^{\tau, \sigma}, l_s^{\tau, \sigma})ds - \pi_{s}^{\tau, \sigma}dW_{s} -\int_{{\bf E}}  l_s^{\tau, \sigma} (e)\tilde{N}(ds,de);
&& X_{\tau \wedge \sigma}^{\tau, \sigma}= I(\tau, \sigma).
 \end{eqnarray*}
At time $S$, the first (resp. second) player chooses a stopping time $\tau$ (resp. $\sigma$) greater than $S$, and looks for maximizing (resp. minimizing) the criterium.

 For each stopping time $S \in {\cal T}_0$, the {\em upper} and {\em lower value functions} at time $S$ are defined respectively by \begin{equation} \label{dessus}
\overline{V}(S):=ess \inf_{\sigma \in \T_S } ess \sup_{\tau \in \T_S} \cal{E}_{S,\tau \wedge \sigma}(I(\tau, \sigma));
\end{equation}
\begin{equation} \label{dessous}
\underline{V}(S):=ess\sup_{\tau\in \T_S } ess \inf_{\sigma \in \T_S} \cal{E}_{S,\tau \wedge \sigma}(I(\tau, \sigma)).
\end{equation}
We clearly have the inequality  $\underline{V}(S) \leq \overline{V}(S)$ a.s.

By definition, we say that there {\em exists a value function} at time $S$ for the generalized Dynkin game if $\overline{V}(S)=\underline{V}(S)$ a.s.

We now introduce the definition of an $S$-saddle point for this game problem.
\begin{definition}
Let $S \in \T_0$. A pair $(\tau^*, \sigma^*) \in \T_S^2$ is called an { \em $S$-saddle point } for the generalized Dynkin game if for each $(\tau, \sigma) \in \T_S^2$ we have
$$ {\cal E}_{S, \tau \wedge \sigma^*}(I(\tau, \sigma^*)) \leq {\cal E}_{S, \tau^* \wedge \sigma^*}(I(\tau^*, \sigma^*)) \leq {\cal E}_{S, \tau^* \wedge \sigma}(I(\tau^*, \sigma))  \quad \text{ a.s. }$$
\end{definition}
We first provide a sufficient condition for the existence of an $S$-saddle point and for the characterization of the common value function as  the solution of the DRBSDE.

\begin{lemma}\label{sufficient}
Suppose that the driver $g$ satisfies Assumption $\eqref{Royer}$. Let $\xi$ and $\zeta$ be RCLL adapted processes in ${\cal S}^2$ such that $\xi_T= \zeta_T$ a.s. and $\xi_{t} \leq \zeta_{t}$, $0 \leq t \leq T$  a.s.\, Suppose that Mokobodski's condition is satisfied.\\
Let $(Y,Z,k(\cdot),A,A')$ be the solution of the DRBSDE (\ref{DRBSDE}).
Let $S$ $\in$ $\T_0$. Let $(\hat{\tau}, \hat{\sigma}) \in \mathcal{T}_S$. Suppose that $(Y_t, \,S\leq t \leq \hat{\tau} )$ is a strong ${\cal E}$-submartingale and that $(Y_t, \,S\leq t \leq\hat{\sigma} )$ is a strong ${\cal E}$-supermartingale with $Y_{\hat{\tau}}=\xi_{\hat{\tau}}$ and $Y_{\hat{\sigma}}=\zeta_{\hat{\sigma}}$ a.s.\,\\
The pair $(\hat{\tau}, \hat{\sigma})$ is then an $S$-saddle point  for the generalized Dynkin game \eqref{dessus}- \eqref{dessous} and $$Y_S=\overline{V}(S)=\underline{V}(S)  \; \text {a.s. }$$
\end{lemma}

\dproof 
Since  the process
$(Y_t, \,S\leq t \leq
\hat{\tau} \wedge \hat{\sigma} )$ is a strong ${\cal E}$-martingale (see Definition \ref{defmart})  
and since $Y_{\hat{\tau}}=\xi_{\hat{\tau}}$ and $Y_{\hat{\sigma}}=\zeta_{\hat{\sigma}}$ a.s.\,,
$$Y_S= {\cal E}_{S,\hat{\tau} \wedge \hat{\sigma}}(Y_{\hat{\tau} \wedge \hat{\sigma}})= {\cal E}_{S,\hat{\tau} \wedge \hat{\sigma}}(\xi_{\hat{\tau} }\textbf{1}_{\hat{\tau} \leq \hat\sigma}+\zeta_{\hat\sigma}\textbf{1}_{\hat\sigma < \hat\tau})=  {\cal E}_{S, \hat\tau \wedge \hat\sigma}(I(\hat\tau ,\hat\sigma))\quad \mbox{a.s.}$$
Let  $\tau \in \T_S$. We want to show that for each $\tau \in \T_S$ 
\begin{equation}\label{quatre}
Y_S \geq {\cal E}_{S,\tau \wedge \hat\sigma}(I (\tau, \hat\sigma)) \quad \text{ a.s. }
\end{equation}
Since the process $(Y_t, \,S\leq t \leq
\tau \wedge \hat\sigma)$ 
is a strong ${\cal E}$-supermartingale, we get
\begin{equation}\label{supermart}
Y_S \geq {\cal E}_{S,\tau \wedge \hat\sigma}(Y_{\tau \wedge \hat\sigma}) \quad \text{ a.s. }
\end{equation}
Since $Y \geq \xi$ and $Y_ { \hat\sigma } = \zeta_{  \hat\sigma}$ a.s.\,,
we also have
$$ Y_{\tau \wedge {\hat\sigma}}= Y_{\tau}\textbf{1}_{\tau \leq \hat\sigma}+Y_{\hat\sigma}\textbf{1}_{\hat\sigma < \tau} \geq \xi_{\tau}\textbf{1}_{\tau \leq \hat\sigma}+\zeta_{\hat\sigma}\textbf{1}_{\hat\sigma < \tau}=I(\tau, \hat\sigma)  \quad \text{ a.s. }$$
By inequality \eqref{supermart} and the monotonicity property of ${\cal E}$, we derive inequality $\eqref{quatre}$. \\
Similarly, one can show that for each $\sigma \in \T_S$, we have:
$$Y_S \leq {\cal E}_{S,\hat\tau \wedge \sigma}(I (\hat\tau, \sigma)) \quad \mbox{a.s.}$$
The pair $(\hat\tau, \hat\sigma)$ is thus an $S$-saddle point and $Y_S=\overline{V}(S)=\underline{V}(S) $ a.s.\,
\fproof


\begin{remark}
When $f$ does not depend on $y,z,k$, from the above proposition, one can derive a well-known sufficient condition of optimality for classical Dynkin game problems ( see e.g. Theorem 2.4 in \cite{ALM} or Proposition 3.1. in \cite{KQC}).
\end{remark}
We now provide an existence result under an additional assumption.
\begin{theorem}[Existence of $S$-saddle points]\label{optimal}
Suppose that $g$ satisfies Assumption $\eqref{Royer}$. Let $\xi$ and $\zeta$ be RCLL adapted processes in ${\cal S}^2$ such that $\xi_{T} = \zeta_{T}$ a.s.\,and $\xi_{t} \leq \zeta_{t}$, $0 \leq t \leq T$  a.s.\, Suppose that Mokobodski's condition is satisfied.\\
Let $(Y,Z,k, A,A')$ be the solution of the DRBSDE (\ref{DRBSDE}). Suppose that $A,A'$ are continuous (which is the case if $ \xi $ and $ -\zeta$ are  l.u.s.c. along stopping times).
For each $S$ $\in$ $\T_0$, consider
$$\tau^*_S:= \inf \{ t \geq S,\,\, Y_t = \xi_t\} ;  \quad \, \sigma^*_S:= \inf \{ t \geq S,\,\, Y_t = \zeta_t\}.$$
$$ \overline\tau_S:= \inf \{ t \geq S,\,\, A_t > A_S\} ;  \quad \overline\sigma_S:= \inf \{ t \geq S,\,\, A'_t > A'_S\}.$$
Then, for each $S \in \T_0$, the pairs of stopping times $(\tau_S^*, \sigma_S^*)$ and $(\overline{\tau}_S, \overline{\sigma}_S)$ are  $S$-saddle points for the generalized Dynkin game and $Y_S=\overline{V}(S)=\underline{V}(S)  \; \text {a.s. }$\\
Moreover, $Y_ { \sigma^*_S } = \zeta_{  \sigma^*_S}$, $Y_ { \tau^*_S } = \xi_{  \tau^*_S}$, $A_{\tau_S^*}=A_S$ and $A'_{\sigma_S^*}=A'_S$ a.s. The same properties hold for $\overline{\tau}_S, \overline{\sigma}_S$.

\end{theorem}

\begin{remark} \label{susu} Note that $\sigma_S^* \leq \overline\sigma_S$ and $\tau_S^* \leq \overline\tau_S$ a.s.\, 
Moreover, by Proposition \ref{compref} in the Appendix,
$(Y_t, \,S\leq t \leq \overline{\tau}_S )$ is a strong ${\cal E}$-submartingale and $(Y_t, \,S\leq t \leq \overline{\sigma}_S )$ is a strong ${\cal E}$-supermartingale. 
\end{remark}
\dproof
Let $S \in \T_0$. Since $Y$ and $\xi$ are right-continuous processes, we have $Y_ { \sigma^*_S } = \xi_{  \sigma^*_S}$ and $Y_ { \tau^*_S } = \xi_{  \tau^*_S}$ a.s.
By definition of $\tau^*_S$, for almost every $\omega$, we have $Y_t(\omega)> \xi_t(\omega)$ for each $t \in$ $[S(\omega), \tau^*_S(\omega)[$.  Hence, since $Y$ is solution of the DRBSDE, the continuous process $A$ is constant on $[S, \tau^*_S]$ a.s.\,because $A$ is continuous. Similarly, the process $A'$ is constant on $[S, \sigma_S^{*}]$ a.s.
By Lemma $\ref{sufficient}$, $( \tau^*_S,  \sigma^*_S)$ is an $S$-saddle point and $Y_S=\overline{V}(S)=\underline{V}(S)  \; \text {a.s. }$\\
It remains to show that $(\overline{\tau}_S, \overline{\sigma}_S)$ is an   $S$-saddle point.
By definition of $\overline{\tau}_S$, $\overline{\sigma}_S$, we have $A_{\overline{\tau}_S}=A_S$ a.s. and $A'_{\overline{\sigma}_S}=A'_S$ a.s. because $A$ and $A'$ are continuous and  $\overline{\tau}_S$, $\overline{\sigma}_S$ are predictable stopping times. Moreover, since the continuous process $A$ increases only on $\{Y_t=\xi_t\}$, we have $Y_{\overline{\tau}_S}=\xi_{\overline{\tau}_S}$ a.s.\, Similarly, $Y_{\overline{\sigma}_S}=\zeta_{\overline{\sigma}_S}$ a.s.\, The result then follows from Lemma \ref{sufficient}.
\fproof

We now see that it is not necessary to have the existence of an $S$-saddle point to ensure the existence of a common value function and its characterization as the solution of a DRBSDE.  
\begin{theorem}[Existence of the value function]\label{caracterisation}

Suppose that $g$ satisfies Assumption $\eqref{Royer}$. Let $\xi$ and $\zeta$ be RCLL adapted processes in ${\cal S}^2$ such that  $\xi_{T}= \zeta_{T}$ a.s. and $\xi_{t}\leq \zeta_{t}$, $0 \leq t \leq T$  a.s.\,
Suppose that Mokobodski's condition is satisfied. 
Let $(Y,Z,k,A,A')$ be the solution of the DRBSDE (\ref{DRBSDE}).\\
There exists a value function for the generalized Dynkin game, and for each stopping time $S$ $\in$ $\T_0$, we have
 \begin{equation}\label{car}
Y_S = \overline{V}(S)=\underline{V}(S) \quad \mbox{a.s.}
\end{equation}

\end{theorem}

\dproof
For each $S$ $\in$ $\T_0$ and for each $\varepsilon >0$, let
$\tau^{\varepsilon}_S$ and $\sigma^{\varepsilon}_S$ be the stopping times defined by
\begin{equation}\label{tauepsilon}
\tau^{\varepsilon}_S := \inf \{ t \geq S,\,\, Y_t \leq \xi_t + \varepsilon\}.
\end{equation}
\begin{equation}\label{sigmaepsilon}
\sigma^{\varepsilon}_S := \inf \{ t \geq S,\,\, Y_t \geq \zeta_t - \varepsilon\}.
\end{equation}
We first show two lemmas.
\begin{lemma}\label{lala}
\begin{itemize}
\item We have
\begin{equation}\label{tau}
 Y_{\tau^{\varepsilon}_S} \leq \xi_{\tau^{\varepsilon}_S} + \varepsilon \quad \mbox{a.s.}
\end{equation}
\begin{equation}\label{alpha}
 Y_{\sigma^{\varepsilon}_S} \geq \zeta_{\sigma^{\varepsilon}_S} - \varepsilon \quad \mbox{a.s.}
\end{equation}

\item
We have $A_{\tau_S^\varepsilon}=A_S$ a.s. and $A'_{\sigma_S^\varepsilon}=A'_S$ a.s.
\end{itemize}
\end{lemma}

\begin{remark} \label{lalabis}
By the second point and Proposition $\ref{compref}$ in the Appendix, the process
$(Y_t, \,S\leq t \leq
\tau^{\varepsilon}_S )$ is a strong ${\cal E}$-submartingale and the process $(Y_t, \,S\leq t \leq \sigma^{\varepsilon}_S )$ is a strong ${\cal E}$-supermartingale.
\end{remark}
\dproof
The first point follows from the definitions of $\tau^{\varepsilon}_S$ and $\sigma^{\varepsilon}_S$ and the right-continuity of $\xi$, $\zeta$ and $Y$. Let us show the second point.
Note that $\tau^{\varepsilon}_S$ $\in$ $\T_S$ and $\sigma^{\varepsilon}_S$ $\in$ $\T_S$. Fix $\varepsilon >0$.
For a.e. $\omega$, if $t \in [S(\omega), \tau^{\varepsilon}_S(\omega)[$, then $Y_t(\omega)> \xi_t (\omega)+ \varepsilon$ and hence $Y_t(\omega)> \xi_t (\omega)$. It follows that almost surely, $A^c$ is constant on $[S, \tau^{\varepsilon}_S]$ and $ A^d$ is constant on $[S, \tau^{\varepsilon}_S[$.
Also, $ Y_{ (\tau^{\varepsilon}_S)^-  } \geq \xi_{ (\tau^{\varepsilon}_S)^-  }+ \varepsilon \,$ a.s.\,\,
Since $\varepsilon >0$, it follows that $ Y_{ (\tau^{\varepsilon}_S)^-  } > \xi_{ (\tau^{\varepsilon}_S)^-  }\quad \mbox{a.s.}$, which implies that $\Delta A^d _ {\tau^{\varepsilon}_S  } =0$ a.s.\,
Hence, almost surely, $A$ is constant on $[S, \tau^{\varepsilon}_S]$.
Similarly, one can show that  $A^{'}$ is a.s.\, constant on $[S, \sigma^{\varepsilon}_S]$.
\fproof

\begin{lemma}\label{eps}
Let $\varepsilon$ $>0$. For all $S$ $\in$ $\T_0$ and $(\tau, \sigma) \in \T_S^2$, we have

\begin{equation}\label{fifi}
{\cal E}_{S, \tau \wedge \sigma^{\varepsilon}_S} ( I(\tau,  \sigma^{\varepsilon}_S))  - K\varepsilon  \,\,\leq \,\,  Y_S \,\,\leq \,\,  {\cal E}_{S, \tau^{\varepsilon}_S \wedge \sigma}( I(\tau^{\varepsilon}_S,  \sigma))  + K\varepsilon \quad \mbox{a.s.}\,,
\end{equation}
where $K$ is a positive constant which only depends on $T$ and the Lipschitz constant $C$ of $f$.

 \end{lemma}

\dproof
Let $\tau \in {\cal T}_S.$ By Remark $\ref{susu}$, the process $(Y_t, \,S\leq t \leq \sigma^{\varepsilon}_S )$ is a strong ${\cal E}$-supermartingale. Hence,
\begin{equation}\label{lala1}
Y_S \geq {\cal E}_{S, \tau \wedge \sigma^{\varepsilon}_S}(Y_{\tau \wedge \sigma^{\varepsilon}_S}) \quad \mbox{a.s.}\
\end{equation}
Since $Y \geq \xi$  and  $Y_{\sigma^{\varepsilon}_S} \geq \zeta_{\sigma^{\varepsilon}_S} - \varepsilon $ a.s. (see Lemma \ref{lala}), we have:
\begin{equation*}
Y_{\tau \wedge \sigma^{\varepsilon}_S} \geq \xi_{\tau} {\bf 1} _{\tau \leq \sigma^{\varepsilon}_S}+ (\zeta_{\sigma^{\varepsilon}_S} -\varepsilon){\bf 1} _{\sigma^{\varepsilon}_S < \tau} \geq I(\tau, \sigma^{\varepsilon}_S) -\varepsilon \quad {\rm a.s.}
\end{equation*}
where the last inequality follows from the definition of $I(\tau, \sigma)$.
Hence, using $\eqref{lala1}$ and the monotonicity property of $\cal E$, we  get
\begin{equation}
Y_{S} \geq  {\cal E}_{S, \tau \wedge \sigma^{\varepsilon}_S}( I(\tau, \sigma^{\varepsilon}_S) -\varepsilon) \quad \mbox{a.s.}
\end{equation}
Now, by the a priori estimates on BSDEs (see Proposition A.4, \cite{16}), we have
\begin{equation*}
 \vert {\cal E}_{S, \tau \wedge \sigma^{\varepsilon}_S}( I(\tau, \sigma^{\varepsilon}_S) -\varepsilon)
- {\cal E}_{S, \tau \wedge \sigma^{\varepsilon}_S}( I(\tau, \sigma^{\varepsilon}_S) ) \vert
\leq K \varepsilon\,\quad \mbox{a.s.}
\end{equation*}
  It follows that
 \begin{equation*}
 Y_{S}\,\, \geq \,\,{\cal E}_{S, \tau \wedge \sigma^{\varepsilon}_S} ( I(\tau,  \sigma^{\varepsilon}_S))  - K\varepsilon   \,\, \quad \mbox{a.s.}\,
\end{equation*}
Similarly, one can show that 
\begin{equation*}
Y_S \,\,\leq \,\,  {\cal E}_{S, \tau^{\varepsilon}_S \wedge \sigma}( I(\tau^{\varepsilon}_S,  \sigma))  + K\varepsilon \quad \mbox{a.s.}\,,
\end{equation*}
which ends the proof of Lemma \ref{eps}.

\fproof

 {\bf End of proof of Theorem \ref{caracterisation}}\\
Using Lemma \ref{eps}, we derive that for each $\varepsilon >0$, 
$$ ess\sup_{\tau \in \T_s}{\cal E}_{S, \tau \wedge \sigma^{\varepsilon}_S}( I(\tau, \sigma^{\varepsilon}_S)) -K \varepsilon \,\,\leq \,\,Y_S \leq ess\inf_{\sigma \in \T_S} {\cal E}_{S, \tau \wedge \sigma^{\varepsilon}_S}( I(\tau, \sigma^{\varepsilon}_S))+K\varepsilon \,\,\text{ a.s.}\,,$$
which implies 
$$\overline{V}(S)- K \varepsilon \,\, \leq \,\, Y_S \,\,\leq \,\, \underline{V}(S)+ K \varepsilon \quad \text{ a.s.}$$
Since 
$\underline{V}(S) \leq \overline{V}(S)$ a.s.\,, we get 
$\underline{V}(S)= Y_S = \overline{V}(S) \text{ a.s.}\,$
The proof of Theorem \ref{caracterisation} is thus complete.
\fproof
\begin{remark}
Inequality $\eqref{fifi}$ shows that $(\tau_S^{\varepsilon}, \sigma_S^{\varepsilon})$ defined by $\eqref{tauepsilon}$  and  $\eqref{sigmaepsilon}$ is an $\varepsilon'$-saddle point at time $S$ with $\varepsilon'=K \varepsilon$.
\end{remark}

\subsection{Generalized mixed game problems}\label{mixed}

We now consider  a generalized mixed game problem when the players have two actions: continuous control and stopping.

Let $(g^{u,v}; (u,v) \in \mathcal{U} \times \mathcal{V})$ be a family of Lipschitz drivers satisfying Assumption \eqref{Royer} .

Let $S \in \mathcal{T}_0.$ For each quadruple $(u,\tau,v,\sigma) \in \mathcal{U} \times \mathcal{T}_S \times \mathcal{V} \times \mathcal{T}_S$, the {\em criterium} at time $S$ is given by $\mathcal{E}_{S, \tau \wedge \sigma}^{u,v}(I(\tau, \sigma)),$ where $\mathcal{E}^{u,v}$ corresponds to the $g^{u,v}$-conditional expectation. The first (resp. second) player chooses a pair $(u, \tau)$ (resp. $(v, \sigma)$) of control and stopping time, and looks for maximizing (resp. minimizing) the criterium.

For each stopping time $S \in {\cal T}_0$, the {\em upper} and {\em lower value functions} at time $S$ are defined respectively by 
\begin{equation} \label{dessusbis}
\overline{V}(S):=ess \inf_{v \in \mathcal{V}, \sigma \in \T_S } ess \sup_{u \in \mathcal{U},\tau \in \T_S} \cal{E}^{u,v}_{S,\tau \wedge \sigma}(I(\tau, \sigma));
\end{equation}

\begin{equation}\label{dessousbis}
\underline{V}(S):=ess\sup_{u \in \mathcal{U},\tau\in \T_S } ess \inf_{v \in \mathcal{V}, \sigma \in \T_S} \cal{E}^{u,v}_{S,\tau \wedge \sigma}(I(\tau, \sigma)).
\end{equation}

We say that there {\em exists a value function} at time $S$ for the game problem if $\overline{V}(S)=\underline{V}(S)$ a.s.
We now introduce the definition of an $S$-saddle point for this game problem.
\begin{definition}
Let $S \in \T_0$. A quadruple $(\overline{u},\overline{\tau}, \overline{v}, \overline{\sigma}  ) \in  \mathcal{U} \times \T_S \times \mathcal{V} \times \T_S$ is called an { \em $S$-saddle point } for the generalized mixed game problem if for each $(u,\tau, v, \sigma  ) \in \mathcal{U} \times \T_S \times \mathcal{V} \times \T_S$ we have
$$ {\cal E}^{{u},\overline{v}}_{S, \tau \wedge  \overline{\sigma}}(I(\tau,  \overline{\sigma})) \leq {\cal E}^{\overline{u},  \overline{v}}_{S, \overline{\tau} \wedge\overline{\sigma}}(I(\overline{\tau} \wedge\overline{\sigma})) \leq {\cal E}^{\overline{u}, v}_{S, \overline{\tau} \wedge \sigma}(I(\overline{\tau}, \sigma))  \quad \text{ a.s. }$$
\end{definition}
We now show that when the obstacles are supposed to be l.u.s.c. along stopping times, there exist some saddle points for the above generalized mixed game problem.
\begin{theorem}\label{generalize}
Let $(g^{u,v}; (u,v) \in \mathcal{U} \times \mathcal{V})$ be a family of Lipschitz drivers satisfying Assumptions $\eqref{Royer}$. Let $\xi$ and $\zeta$ be RCLL  adapted processes in ${\cal S}^2$ and l.u.s.c. along stopping times, such that $\xi_{T} = \zeta_{T}$ a.s. and $\xi_{t} \leq \zeta_{t}$, $0 \leq t \leq T$  a.s.\, Suppose that Mokobodski's condition is satisfied and that there exist controls $\overline{u} \in \mathcal{U}$ and $\overline{v} \in \mathcal{V}$ such that for each $(u,v) \in \mathcal{U} \times \mathcal{V},$
\begin{equation}\label{Issac}
g^{ u,\overline{v}}(t, Y_t,Z_t, k_t) \leq g^{\overline{u},\overline{v}}(t,Y_t,Z_t,k_t) \leq g^{ \overline{u},v}(t,Y_t,Z_t,k_t) \quad dt \otimes dP \text{ a.s. },
\end{equation}
where $(Y,Z,k,A,A')$ corresponds to the solution of the DRBSDE (\ref{DRBSDE}) associated with driver $g^{\overline{u},\overline{v}}$.
Consider the stopping times
\begin{equation*}
\tau_S^*:= \inf \{t \geq S: Y_t = \xi_t \} \quad; \quad \sigma_S^*:= \inf \{t \geq S: Y_t = \zeta_t \}.
\end{equation*}
The quadruple $(\overline{u},\tau_S^*, \overline{v}, \sigma_S^*  )$ is then an S-saddle point for the generalized mixed game problem \eqref{dessusbis}-\eqref{dessousbis}, and  we have $Y_S=\underline{V}(S)=\overline{V}(S)$ a.s.
\end{theorem}

\dproof
By the last assertion of Theoreom $\ref{optimal}$, the process
$(Y_t, \,S\leq t \leq
\tau_S^* \wedge\sigma_S^*)$ is a strong ${\cal E}^{\overline{u},\overline{v}}$-martingale and  $Y_ { \tau_S^*} = \zeta_{  \tau_S^*}$, 
$Y_ { \sigma_S^*} = \zeta_{  \sigma_S^*} $ a.s.\,, which implies
$$Y_S= {\cal E}^{\overline{u},\overline{v}}_{S,\tau_S^* \wedge \sigma_S^*}(Y_{\tau_S^* \wedge \sigma_S^*})= {\cal E}^{\overline{u},\overline{v}}_{S,\tau_S^* \wedge \sigma_S^*}(\xi_{\tau_S^*}\textbf{1}_{\tau_S^*\leq \sigma_S^*}+\zeta_{\sigma_S^*}\textbf{1}_{\sigma_S^* < \tau_S^*})=  {\cal E}^{\overline{u},\overline{v}}_{S, \tau_S^* \wedge \sigma_S^*}(I(\tau_S^* ,\sigma_S^*))\quad \mbox{a.s.}$$
Let  $\tau \in \T_S$. 
Since $Y\geq \xi$ and $Y_ { \sigma_S^*} = \zeta_{  \sigma_S^*}$ a.s.\,,
we have
$$ Y_{\tau \wedge {\sigma_S^*}}= Y_{\tau}\textbf{1}_{\tau \leq \sigma_S^*}+Y_{\sigma_S^*}\textbf{1}_{\sigma_S^* < \tau} \geq \xi_{\tau}\textbf{1}_{\tau \leq \sigma_S^*}+\zeta_{\sigma_S^*}\textbf{1}_{\sigma_S^* < \tau}=I(\tau, \sigma_S^*)\quad \mbox{a.s.}$$
Moreover, by Theorem \ref{optimal}, $A'_{\sigma_S^*}=A'_s$ a.s., which implies that:
\begin{equation*}
-dY_t = g^{\overline{u}, \overline{v}}(t,Y_t,  Z_t, k_t)dt +dA_t -Z_t  dW_t -\int_{{\bf E}}  k_t(e) \tilde{N}(dt,de); \quad S\leq t \leq \sigma_S^*, \quad dt\otimes dP \text{ a.s.}
\end{equation*}
Hence, $(Y_t)_{S \leq t \leq \tau \wedge \sigma_S^*}$ is the solution of the BSDE associated with generalized driver $``g^{\overline{u},\overline{v}}(\cdot)dt+dA_t"$ and terminal condition $Y_{\tau \wedge \sigma_S^*}$ . By using Assumption $\eqref{Issac}$, the inequality $Y_{\tau \wedge \sigma_S^*} \geq I(\tau, \sigma_S^*)$ and the comparison theorem for BSDEs with jumps, we obtain that  for each $u \in \mathcal{U}$:
\begin{equation*}
Y_S\geq \mathcal{E}^{u,\overline{v}}_{S, \tau \wedge \sigma_S^* }(I(\tau,\sigma_S^*)) \quad \mbox{a.s.}
\end{equation*}
Similarly, one can show that for each $v \in \mathcal{V }, \sigma \in \T_S$, we have:
$$Y _S \leq {\cal E} ^{\overline{u},v}_{S,\tau_S^* \wedge \sigma}(I (\tau_S^*, \sigma))\quad \mbox{a.s.}$$
The quadruple $(\overline{u},\tau_S^*, \overline{v}, \sigma_S^*  )$ is thus an $S$-saddle point and $Y_S=\overline{V}(S)=\underline{V}(S) $ a.s.\,
\fproof

Under less restricted assumptions on the obstacles, we show that there exist a value function for the above game problem which can be characterized as the solution of a DRBSDE.
\begin{theorem}[Existence of the value function]\label{caracterisation2}

Let $(g^{u,v}; (u,v) \in \mathcal{U} \times \mathcal{V})$ be a family of drivers satisfying Assumptions $\eqref{Royer}$ and which are uniformly Lipschitz with common Lipchitz constant $C$. Let $\xi$ and $\zeta$ be RCLL adapted processes in ${\cal S}^2$ such that  $\xi_{T}= \zeta_{T}$ a.s.\,and $\xi_{t}\leq \zeta_{t}$, $0 \leq t \leq T$  a.s.\,
Suppose that Mokobodski's condition is satisfied and that there exist controls $\overline{u} \in \mathcal{U}$ and $\overline{v} \in \mathcal{V}$ such that for each $u \in \mathcal{U}, v \in \mathcal{V}$:
\begin{equation}\label{Issac1}
g^{ u,\overline{v}}(t, Y_t,Z_t, k_t) \leq g^{\overline{u},\overline{v}}(t,Y_t,Z_t,k_t) \leq g^{ \overline{u},v}(t,Y_t,Z_t,k_t), \quad dt \otimes dP \text{ a.s. }
\end{equation}associated with
where $(Y,Z,k,A,A')$ corresponds to the solution of the DRBSDE (\ref{DRBSDE}) associated with driver $g^{\overline{u},\overline{v}}$.\\
Then, there exists a value function for the generalized mixed game problem \eqref{dessusbis}-\eqref{dessousbis}, and for each stopping time $S$ $\in$ $\T_0$, we have
 \begin{equation*}
Y_S = \overline{V}(S)=\underline{V}(S) \quad \mbox{a.s.}
\end{equation*}

\end{theorem}

\dproof
For each $S$ $\in$ $\T_0$ and for each $\varepsilon >0$, let
$\tau^{\varepsilon}_S$ and $\sigma^{\varepsilon}_S$ be the stopping times defined by
$$\tau^{\varepsilon}_S := \inf \{ t \geq S,\,\, Y_t \leq \xi_t + \varepsilon\}; \\ \quad
\sigma^{\varepsilon}_S := \inf \{ t \geq S,\,\, Y_t \geq \zeta_t - \varepsilon\}.$$
Let $\tau \in {\cal T}_S$. Since $Y \geq \xi$  and  $Y_{\sigma^{\varepsilon}_S} \geq \zeta_{\sigma^{\varepsilon}_S} - \varepsilon $ a.s. ( see Lemma \ref{lala}), we have:
\begin{equation*}
Y_{\tau \wedge \sigma^{\varepsilon}_S} \geq \xi_{\tau} {\bf 1} _{\tau \leq \sigma^{\varepsilon}_S}+ (\zeta_{\sigma^{\varepsilon}_S} -\varepsilon){\bf 1} _{\sigma^{\varepsilon}_S < \tau}\geq  I(\tau, \sigma^{\varepsilon}_S) -\varepsilon \quad {\rm a.s.}
\end{equation*}
By Lemma $\ref{lala}$, $A'_{\sigma_S^{\varepsilon}}=A'_S$ a.s. which implies that:
\begin{equation*}
-dY_t = g^{\overline{u}, \overline{v}}(t,Y_t,  Z_t, k_t)dt +dA_t -Z_t  dW_t -\int_{{\bf E}}  k_t(e) \tilde{N}(dt,de), \\ \quad S \leq t \leq \sigma_S^{\varepsilon}, \quad dt \otimes dP \text{ a.s.}
\end{equation*}
Hence, $(Y_t)_{S \leq t \leq \tau \wedge \sigma^{\varepsilon}}$ is the solution of the BSDE associated with generalized driver $``f(\cdot)dt+dA_t"$ and terminal condition $Y_{\tau \wedge \sigma^\varepsilon}$. By using Assumption $\eqref{Issac1}$, the inequality $Y_{\tau \wedge \sigma^\varepsilon} \geq I(\tau, \sigma^\varepsilon)-\varepsilon$ and the comparison theorem for BSDEs with jumps, we obtain
\begin{equation*}
Y_S\geq \mathcal{E}^{u,\overline{v}}_S(I(\tau,\sigma^{\varepsilon})-\varepsilon)\geq \mathcal{E}^{u,\overline{v}}_S(I(\tau,\sigma^{\varepsilon}))-K\varepsilon   \quad \text{  a.s. },
\end{equation*}
where the second inequality follows from the a priori estimates for BSDEs with jumps.
Here, the constant $K$ only depends on $T$ and $C$, the common Lipschitz constant.
Consequently, we get
\begin{equation*}
Y_S \geq ess \inf_{v \in \mathcal{V}, \sigma \in \T_S } ess \sup_{u \in \mathcal{U},\tau \in \T_S} \cal{E}^{u,v}_{S,\tau \wedge \sigma}(I(\tau, \sigma))-K \varepsilon  \quad  \text{ a.s.}
\end{equation*}
Similarly, one can show that for each $\varepsilon >0$,
\begin{equation*}\label{eq2}
Y_S\leq ess \sup_{u \in \mathcal{U},\tau \in \T_S}  ess \inf_{v \in \mathcal{V}, \sigma \in \T_S }\cal{E}^{u,v}_{S,\tau \wedge \sigma}(I(\tau, \sigma))+K \varepsilon \quad  \text{ a.s.}
\end{equation*}
Hence, $\overline{V}(S) \leq \underline{V}(S)$  a.s.\, Since $\underline{V}(S) \leq \overline{V}(S)$ a.s., 
the equality follows.
\fproof


\begin{remark}
Note that Theorem \ref{caracterisation2} still holds if $g^{\overline{u},\overline{v}}$ is replaced by any Lipschitz driver $g$ which satisfies \eqref{Issac1}.
\end{remark}
\paragraph{Application:}

\bigskip

\quad  Let $U,V$ be compact polish spaces.\\
We are given a map $F:[0,T] \times \Omega \times U \times V \times \R^2 \times L_{\nu}^2 \rightarrow \R,$ $(t, \omega,u,v,y,z,k) \mapsto F(t, \omega,u,v,y,z,k)$, supposed to be measurable with respect to $\mathcal{P} \otimes \mathcal{B}(U) \otimes \mathcal{B}(V) \otimes \mathcal{B}(\R^2) \otimes \mathcal{B}(L_{\nu}^2),$ continuous, concave (resp. convex) with respect to $u$ ( resp. $v$), and uniformly Lipchitz with respect to $(y,z,k)$. Suppose that $F(t, \omega, u,v,0,0,0)$ is uniformly bounded.\\
\quad Let $\mathcal{U}$ (resp. $\mathcal{V}$) be the set of predictable processes valued in $U$ ( resp. $V$). For each $(u,v) \in \mathcal{U} \times \mathcal{V}$, let $g^{u,v}$ be the driver defined by $$g^{u,v}(t,\omega,y,z,k):=
F(t, \omega, u_t(\omega),v_t(\omega),y,z,k).$$ 
 Let $\xi$ and $\zeta$ be RCLL adapted processes in ${\cal S}^2$ such that  $\xi_{T}= \zeta_{T}$ a.s.\,and $\xi_{t}\leq \zeta_{t}$, $0 \leq t \leq T$  a.s.\, Suppose that Mokobodski's condition is satisfied.\\
Let us consider the associated generalized mixed game problem. Define for each $(t, \omega,y,z,k)$ the map 
\begin{equation} \label{ff}
g(t,\omega,y,z,k)=\sup_{u \in U} \inf_{v \in V} F(t, \omega,u,v,y,z,k).
\end{equation}
 Since $U$ and $V$ are polish spaces, there exist some dense countable subsets $\overline{U}$ (resp. $\overline{V}$) of $U$ (resp. $V$). Since $F$ is continuous with respect to $u,v$, the sup and the inf can be taken over $\overline{U}$ (resp. $\overline{V}$). Hence, $g$ is a Lipchitz driver.\\
\quad Let $(Y,Z,k, A,A') \in \mathcal{S}^2 \times \mathbb{H}^2 \times \mathbb{H}^2_{\nu} \times (\mathcal{A}^2)^2$ be the solution of the DRBSDE associated with driver $g$ and obstacles $\xi$ and $\zeta$. By classical convex analysis, for each $(t, \omega)$  there exist $(u^{*},v^{*}) \in (U,V)$ such that
\begin{eqnarray}\label{mes}
F(t, \omega, u,v^*,Y_{t^-}(\omega),Z_t(\omega),k_t(\omega))& \leq &F(t, \omega, u^*,v^*,Y_{t^-}(\omega),Z_t(\omega),k_t(\omega)) \\ 
& \leq& F(t, \omega, u^*,v,Y_{t^-}(\omega),Z_t(\omega),k_t(\omega)), \text{ } \forall (u,v) \in U \times V ;\nonumber\\
g(t, \omega,Y_{t^-}(\omega),Z_t(\omega),k_t(\omega)))&=&F(t,\omega, u^*,v^*,Y_{t^-}(\omega),Z_t(\omega),k_t(\omega))\nonumber
\end{eqnarray}
Since the set of all $(t, \omega, u^*,v^*) \in [0,T] \times \Omega \times U \times V$ satisfying conditions \eqref{mes} belongs to $\mathcal{P} \times \mathcal{B}(U) \times \mathcal{B}(V)$, by applying 
 the section theorem (see Section 81 in the Appendix of Ch. III in \cite{DM1}), we get that there exists a pair of predictable process $(u^*, v^*) \in \mathcal{U} \times \mathcal{V}$ such that $dt \otimes dP$ a.s., for all $(u,v) \in \mathcal{U} \times \mathcal{V}$ we have $dt \otimes dP$ a.s.:
$$F(t, u_t,v_t^*,Y_{t},Z_t, k_t) \leq F(t,u_t^*,v_t^*,Y_{t},Z_t,k_t) \leq F(t, u_t^*,v_t,Y_{t},Z_t,k_t) $$
and $g(t, Y_{t}, Z_t, k_t)=F(t,u_t^*,v_t^*,Y_{t},Z_t,k_t)$.
Hence, Assumption $\eqref{Issac}$ is satisfied. 
By applying  Theorems \ref{caracterisation2} and \ref{generalize}, we derive the following result: 
\begin{proposition}\label{Hamadene}
There exists a value function for the above generalized mixed game problem (associated with the map $F(t, u,v,y,z,k)$). 
Let $Y$ be the solution of the DRBSDE associated with obstacles $\xi$, $\zeta$ and the driver $g$ defined by \eqref{ff}.\\
For each stopping time $S$ $\in$ $\T_0$, we have $Y_S = \overline{V}(S)=\underline{V}(S)$ a.s.\,\\
Suppose that $\xi$ and $\zeta$ are l.u.s.c. along stopping times. 
Consider the stopping times
\begin{equation*}
\tau_S^*:= \inf \{t \geq S: Y_t = \xi_t \} \quad; \quad \sigma_S^*:= \inf \{t \geq S: Y_t = \zeta_t \}.
\end{equation*}
The quadruple $(u^*,\tau_S^*, v^*, \sigma_S^*  )$ is then an $S$-saddle point for this mixed game problem. 
\end{proposition}
\bigskip

\noindent We give now an example for which the above proposition can be applied.\\
\textbf{Example:}
Let us now consider the particular case when $F$ takes the following form:\\ $F(t, \omega, u,v,y,z,k)=\beta(t,\omega,u,v)y+<\gamma(t,\omega,u,v,\cdot),k>_{\nu}+c(t,\omega,u,v),$ with $\beta, \gamma, c$ bounded.  
By classical results on linear BSDEs (see \cite{16}), the criterium can be written $$\mathcal{E}^{u,v}_{S, \tau \wedge \sigma}(I(\tau, \sigma))=E_{Q^{u,v}}\left[\int_S^{\tau \wedge \sigma}c(t,u_t,v_t)dt+I(\tau, \sigma) | \mathcal{F}_S \right],$$ with $Q^{u,v}$  the probability measure which admits $Z^{u,v}_T$ as density with respect to $P$, where $(Z_t^{u,v})$ is the solution of the following SDE:
$$
dZ_t^{u,v}=Z_t^{u,v}[\beta(t,u_t,v_t)dW_t+\int_{{\bf E}} \gamma(t,u_t,v_t,e)\Tilde{N}(dt,de)]; \quad Z^{u,v}_0=1.
$$
The process $c(t,u_t,v_t)$ can be interpreted as an instantaneous reward associated with controls $u,v$. This linear model takes into account some ambiguity on the model via the probability measures $Q^{u,v}$ as well as some ambiguity on the instantaneous reward.
This case corresponds to the classical mixed game problems studied in \cite{BL} and \cite{16}. The above proofs provide some alternative short proofs of their results.


\section{Comparison theorems for DRBSDEs with jumps and a priori estimates}\label{sec5}

\subsection{Comparison theorems}

\begin{Theorem}[Comparison theorem for DRBSDEs.]\label{thmcomprbsde}
 Let $\xi^{1}$, $\xi^{2}$, $\zeta^{1}$, $\zeta^{2}$ be  processes in ${\cal S}^2$ such that $\xi_T^i = \zeta_T^i$ a.s.\, and $\xi_t^i \leq \zeta_t^i$, $0 \leq t \leq T$ a.s.  for $i=1,2$. Suppose that for $i=1,2$, $\xi^{i}, \zeta^{i}$ satisfies Mokobodski's condition. Let $g^1$and $g^2$
 be  Lipschitz drivers   satisfying Assumption~\eqref{Royer}.\\
Suppose  that
\begin{itemize}
\item   $\xi_t ^{2}\le \xi_t ^{1}$ and $\zeta_t ^{2}\le \zeta_t ^{1}$, $0\leq t \leq T$ a.s.

\item  $g^{2}(t,y,z,k) \le g^{1}(t,y,z,k),\,\, \text{ for all  }  (y,z,k) \in \R^2 \times {\cal L}^2_\nu; \;\;  dP\otimes dt-a.s.\,$
\end{itemize}
Let $(Y^{i}, Z^{i}, k^i, A^i, A^{' i})$  be
the solution of the DRBSDE associated with $(\xi ^{i},\zeta ^{i},g^{i})$ , $i=1,2$. Then,
$$Y_{t}^{2}\le
Y_{t}^{1}, \,\,0\leq t \leq T \quad \text{ a.s. } $$
\end{Theorem}

\begin{Remark}
Note that a comparison theorem has been provided in \cite{CM} in the case of jumps under stronger assumptions.Their proof is different and based on Ito's calculus.
\end{Remark}
\dproof We give a short proof  based on the characterization of solutions of DRBSDEs (Theorem \ref{caracterisation}) via generalized Dynkin games.
 Let $t$ $\in$ $[0,T]$. For each $\tau, \sigma \in T_{t}$, let us denote by ${\cal E}^{i}_{., \tau \wedge \sigma}( I^i(\tau,\sigma))$ the unique solution 
 of
the BSDE associated with driver $g^{i}$, terminal time $\tau \wedge \sigma$ and terminal condition $I^i(\tau,\sigma):=\xi_{\tau}^{i} \textbf{1}_{\tau \leq \sigma}+\zeta_{\sigma}^{i} \textbf{1}_{\sigma < \tau}$ for $i=1,2$.
Since $g^2 \leq g^1$, and $I^2(\tau, \sigma) \leq I^1(\tau, \sigma)$, by the comparison theorem 
for BSDEs, the following inequality
$${\cal E}^{2}_{t, \tau \wedge \sigma}( I^2(\tau,\sigma)) \leq {\cal E}^{1}_{t, \tau \wedge \sigma}( I^1(\tau,\sigma)) \,\,a.s.\,$$
holds for each $\tau$, $\sigma$ in $\T_t$.
Hence, by taking the essential supremum over $\tau$ in $\T_t$ and the essential infimum over $\sigma$ in $\T_t$, and by using Theorem \ref{caracterisation}, we get
$$Y^{2}_{t}= {\rm ess} \inf_{\sigma \in \T_t}{\rm ess} \sup_{\tau \in \T_t}{\cal E}^{2}_{t, \tau \wedge \sigma}( I^2(\tau,\sigma)) \leq {\rm ess} \inf_{\sigma \in \T_t}{\rm ess} \sup_{\tau \in \T_t}{\cal E}^{1}_{t, \tau \wedge \sigma}( I^1(\tau,\sigma))=Y^{1}_{t}\,\,a.s.  $$
\fproof

We now provide a strict comparison theorem, which had not been given in the literature even in the Brownian case. The first assertion addresses the particular case when the non decreasing processes are continuous and the second one deals with the general case.
\begin{Theorem}[Strict comparison.]\label{sctun}
Suppose that  the assumptions of  Theorem  \ref{thmcomprbsde} hold and that the driver $g^1$
satisfies Assumption \ref{Royer} with $\gamma_t >-1$ in \eqref{robis}.
 Let $S$ in $\T_0$ and suppose that
$Y^1_{S} = Y^2_{S}$ a.s.
\begin{enumerate}
\item
Suppose that $A^i,A^{' i}$, $i=1,2$ are continuous.
For $i=1,2$, let \\
 $\overline \tau_i = \overline \tau_{i,S}:= \inf \{s \geq S ; \, A^i_s > A^i_S  \}$ and $\overline \sigma_i = \overline \sigma_{i,S}:= \inf \{s \geq S ; \, A^{' i}_s > A^{' i}_S \}$. Then
$$Y^1_{t} =Y^2_{t} ,\;\;  S \leq t \leq  \overline \tau_1 \wedge \overline \tau_2 \wedge \overline \sigma_1 \wedge \overline \sigma_2 \; \text{ a.s}. $$
and
\begin{equation}\label{autre}
g^2(t,Y^2_t,Z^2_t,k^2_t) \, =\,  g^{1}(t,Y^2_t,Z^2_t,k^2_t) \, \, \,\,\,  \;\;   S \leq t \leq  \overline \tau_1 \wedge \overline \tau_2 \wedge \overline \sigma_1 \wedge \overline \sigma_2, \,\,
dP\otimes dt-a.s.
\end{equation}
\item
Consider the case when $A^i,A^{' i}$, $i=1,2$ are not necessarily continuous. For $i = 1,2$, define for each $\varepsilon >0$,
$$\tau^{\varepsilon}_i:= \inf \{ t \geq S,\,\, Y^i_t \leq \xi^i_t + \varepsilon \} \; ; \;\; \sigma^{\varepsilon}_i:= \inf \{ t \geq S,\,\, Y^i_t \leq \zeta^i_t - \varepsilon \} .\;\;
$$
Setting $\tilde \tau_i:= \lim_{\varepsilon \downarrow 0} \uparrow \tau^{\varepsilon}_i$ and  $\tilde \sigma_i:= \lim_{\varepsilon \downarrow 0} \uparrow \sigma^{\varepsilon}_i , $ we have
\begin{equation}\label{unoss}
Y^1_{t} =Y^2_{t} , \;\;  S \leq t <  \tilde \tau_1 \wedge \tilde \tau_2 \wedge \tilde \sigma_1 \wedge \tilde \sigma _2 . \quad \text{a.s.}
\end{equation}
 Moreover,  equality \eqref{autre} holds on $[S, \tilde \tau_1 \wedge \tilde \tau_2 \wedge \tilde \sigma_1 \wedge \tilde \sigma _2 ]$.
\end{enumerate}
\end{Theorem}

%
\dproof We adopt the same notation as in the proof of the comparison theorem.\\Suppose first that  $A^i,A^{' i}$, $i=1,2$ are continuous.
By Theorem \ref{optimal}, for $i=1,2$, $(\overline \tau_i, \overline \sigma_i) $ is a saddle point for the game problem associated with $g= g^i$, $\xi=\xi^i$ and $\zeta=\zeta^i$. By Remark \ref{susu}, $(Y^i_t, {S \leq t \leq \overline \tau_i \wedge  \overline \sigma_i})$ is an ${\cal E}^{i}$ martingale. Hence we  have $$Y^i_t = {\cal E}^{i}_{t,\overline \tau_i \wedge \overline \sigma_i}(I(\overline \tau_i ,\overline \sigma_i)), \,\,\, S \leq t \leq \overline \tau_i \wedge \overline \sigma_i \,\,\,{\rm a.s.}$$
Setting $\overline \theta=\overline \tau_1 \wedge \overline \tau_2 \wedge \overline \sigma_1 \wedge \overline \sigma_2 $,
we thus have $$Y^i_t = {\cal E}^{i}_{t,\overline \theta}(Y^{i}_{\overline \theta}), \,\,\, S \leq t \leq \overline \theta \,\,\,{\rm a.s.} \text{ for }i=1,2.$$
By hypothesis, $Y^1_S = Y^2_S$ a.s. Now, we apply the strict comparison theorem for non reflected BSDEs with jumps (see \cite{16}, Th 4.4) for terminal time $\overline \theta$. Hence, we get
$Y^1_{t} =Y^2_{t} ,\;\;  S \leq t \leq \overline \theta \;$ a.s.\,, as well as equality (\ref{autre}), which provides the desired result.

Consider now the general case. \\
Let $\varepsilon $ $>0$.
By  Remark \ref{lalabis}, $(Y^i_t, {S \leq t \leq \tau^{\varepsilon}_i \wedge  \sigma^{\varepsilon}_i  })$ is an ${\cal E}^{i}$ martingale. Hence we  have $$Y^i_t = {\cal E}^{i}_{t,\tau_i^{\varepsilon} \wedge \sigma_i^{\varepsilon}}(I(\tau_i^{\varepsilon} ,\sigma_i^{\varepsilon})), \,\,\, S \leq t \leq \tau^{\varepsilon}_i \wedge \sigma^{\varepsilon}_i \,\,\,{\rm a.s.}$$
By the same arguments as above with $\tau^{*}_1$,$\tau^{*}_2$ and $\sigma^{*}_1$,$\sigma^{*}_2$ replaced by $\tau^{\varepsilon}_1$,$\tau^{\varepsilon}_2$ and $\sigma^{\varepsilon}_1$,$\sigma^{\varepsilon}_2$ respectively, we derive
$Y^1_{t} =Y^2_{t} ,\;\;  S \leq t \leq \tau^{\varepsilon}_1 \wedge \tau_2^{\varepsilon} \wedge \sigma^{\varepsilon}_1 \wedge \sigma_2^{\varepsilon}$  a.s.\,, and equality (\ref{autre}) holds on $[S, \tau^{\varepsilon}_1 \wedge\tau_2^{\varepsilon} \wedge \sigma^{\varepsilon}_1 \wedge \sigma_2^{\varepsilon}]$,  $dt \otimes dP$-a.s.\,
By letting $\varepsilon$ tend to $0$, we obtain the desired result.
\fproof

\bigskip

\quad We now give an application of the above comparison theorem to a control game problem for DRBSDEs.

\begin{proposition}[Control game problem for DRBSDEs]

Suppose that the assumptions of Th. $\ref{caracterisation2}$ hold. For each $(u,v) \in \mathcal{U} \times \mathcal{V}$, let $Y^{u,v}$ be the solution of the DRBSDE $\eqref{DRBSDE}$ associated with driver $g^{u,v}$. Then, for each $S \in \mathcal{T}_0$, $Y_S^{u,\overline{v}} \leq Y_S^{\overline{u}, \overline{v}} \leq Y_S^{\overline{u}, v}$ a.s. 
\end{proposition}

\dproof
By using Assumption \eqref{Issac}   and by applying the comparison theorem for DRBSDEs (Th. \ref{thmcomprbsde}), we get that for each $u \in \mathcal{U}$, $Y_S^{{u},\overline{v}} \leq Y_S^{\overline{u}, \overline{v}}$ a.s.\,Similarly, for all $v \in \mathcal{V}$, we have $Y_S^{\overline{u}, \overline{v}} \leq Y_S^{\overline{u},v}$ a.s.
\fproof

\begin{remark}
From this result, it follows that, under the Assumption \eqref{Issac}, the value function of the above control game problem for DRBSDEs coincides with the one associated with the generalized mixed game problem studied in Section \ref{mixed}.
\end{remark}


\subsection{Some new estimates}
Using the links between generalized Dynkin games and DRBSDEs (see Theorem \ref{caracterisation}), we prove the following estimates.

\begin{proposition} \label{oubli}

Let $\xi^{1}, \xi^{2}, \zeta^{1}, \zeta^{2} \in \mathcal{S}^2$ such that $\xi_T^i = \zeta_T^i$ a.s.\, and $\xi_t^i\leq \zeta_t^i$, $0 \leq t \leq T$ a.s.\, Suppose that for $i=1,2$, $\xi^{i}$ and $\zeta^{i}$ satisfy Mokobodski's condition.
Let $g^1, g^2$ be Lipschitz drivers satisfying Assumption \ref{Royer} with Lipschitz constant $C>0$. For $i=1,2$, let $Y^i$ be the solution of the DRBSDE associated with driver $g^i$, terminal time $T$ and barriers $\xi^i$, $\zeta^i$. For $s \in [0,T]$, let $\overline{Y}:=Y^1-Y^2$, $\overline{\xi}:=\xi^1-\xi^2$, $\overline{\zeta}=\zeta^1-\zeta^2$ and $\overline{g}_s:= \sup_{y,z,k}|g^1(s,y,z,k)-g^2(s,y,z,k)|$. Let $\eta, \beta >0$ be such that $\beta\geq \dfrac{3}{\eta}+2C$ and $\eta \leq \dfrac{1}{C^2}$. Then for each $t$, we have:
\begin{equation}\label{eqA.1}
 \overline{Y}_t^2  \leq e^{\beta (T-t)}E[\sup_{s \geq t} \overline{\xi_s}^2+\sup_{s \geq t} \overline{\zeta_s}^2| \mathcal{F}_t]+ \eta E[\int_t^Te^{\beta (s-t)}{\overline{g}_s^2}ds|\mathcal{F}_t] \textsc{ } a.s.
\end{equation}

\end{proposition}

\begin{remark}
The constants  $\eta$ and $\beta$ are universal, i.e. they do not depend on $T$, $\xi^1, \xi^2, g^1, g^2$. Note that in the previous literature, there does not exist any result providing estimates on DRBSDEs, even in the Brownian case.
\end{remark}

\noindent \textbf{Proof.}
For $i=1,2$ and for each $\tau, \sigma \in \tau_{0}$, let $(X^{i,\tau, \sigma}$, $\pi^{i,\tau, \sigma}, l^{i, \tau, \sigma})$  be the solution of the BSDE associated with driver $g^i$, terminal time $\tau \wedge \sigma$ and terminal condition $I^i(\tau, \sigma)$, where $I^i(\tau, \sigma)=\xi_{\tau}^i \textbf{1}_{ \tau \leq \sigma}+\zeta_{\sigma}^i \textbf{1}_{\sigma<\tau}$.
Set  $\overline{X}^{\tau, \sigma}:=X^{1,\tau, \sigma}-X^{2,\tau, \sigma}$ and $\overline{I}^{\tau, \sigma}:=I^1(\tau, \sigma)-I^2(\tau, \sigma)=\overline{\xi}_{\tau} \textbf{1}_{\tau \leq \sigma}+\overline{\zeta}_{\sigma}\textbf{1}_{\sigma<\tau}$.\\
By a priori estimate on BSDEs (see Proposition $A.4$ in \cite{17}), we have a.s.:
\begin{equation}\label{A.2}
 (\overline{X}_t^{\tau, \sigma})^2  \leq e^{\beta (T-t)} E[\overline{I}(\tau, \sigma)^2 \mid \mathcal{F}_t]+ \eta  E[\int_t^T e^{\beta (s-t)}{[(g^1-g^2)(s, X_s^{2,\tau, \sigma},\pi_s^{2,\tau, \sigma}, l_s^{2,\tau, \sigma}) ]^2ds} \mid \mathcal{F}_t]
\end{equation}
from which we derive that
\begin{equation}\label{A.3}
 (\overline{X}_t^{\tau, \sigma})^2  \leq e^{\beta (T-t)} E[\sup_{s \geq t}\overline{\xi}_{s}^2+\sup_{s \geq t}\overline{\zeta}_{s}^2|\mathcal{F}_t]+ \eta E[\int_t^Te^{\beta (s-t)}{\overline{g}_s^2 ds}| \mathcal{F}_t] \quad a.s.
\end{equation}
Now, by using inequality $\eqref{fifi}$, we obtain that for each $\varepsilon >0$ and for all stopping times $\tau, \sigma$,
$$
Y_t^1-Y_t^2 \leq X_t^{1, \tau^{\epsilon, 1}, \sigma}-X_t^{2, \tau, \sigma^{\epsilon, 2}}+2K\epsilon.
$$
Applying this inequality  to  $\tau=\tau^{\epsilon,1},\sigma=\sigma^{\epsilon, 2}$ we get  
\begin{align}\label{ineq2}
Y_t^1-Y_t^2 \leq X_t^{1, \tau^{\epsilon, 1}, \sigma^{\epsilon, 2}}-X_t^{2, \tau^{\epsilon, 1}, \sigma^{\epsilon, 2}}+2K\epsilon
\leq |X_t^{1, \tau^{\epsilon, 1}, \sigma^{\epsilon, 2}}-X_t^{2, \tau^{\epsilon, 1}, \sigma^{\epsilon, 2}}|+2K\epsilon.
\end{align}
By $\eqref{A.3}$ and $\eqref{ineq2}$, we have:
\begin{align*}
Y_t^1-Y_t^2 \leq \sqrt{ e^{\beta (T-t)} E[\sup_{s \geq t} \overline{\xi_s}^2+\sup_{s \geq t} \overline{\zeta_s}^2| \mathcal{F}_t]+ \eta E[\int_t^Te^{\beta (s-t)}{\overline{g}_s^2}ds|\mathcal{F}_t]}+2K\epsilon.
\end{align*}
By symmetry, the last inequality is also verified by $Y_t^2-Y_t^1$.
The result follows.
\fproof

We also provide the following estimate on the common value function $Y$ of our generalized Dynkin game problem (\eqref{dessus} and \eqref{dessous}) (or equivalently the solution of the DRBSDE associated with driver $g$).

\begin{proposition} \label{A.4}

For each $t$, we have:
\begin{equation}
 Y_t^2  \leq e^{\beta (T-t)}E[\sup_{s \geq t} {\xi_s}^2+\sup_{s \geq t} {\zeta_s}^2| \mathcal{F}_t]+ \eta E[\int_t^Te^{\beta (s-t)}{g(s,0,0,0)^2}ds|\mathcal{F}_t] \textsc{ } a.s.
\end{equation}

\end{proposition}

\noindent\textbf{Proof.}
Let $X_t^{\tau, \sigma}$ be the solution of the BSDE associated with driver $g$, terminal time $\tau \wedge \sigma$ and terminal condition $I(\tau, \sigma)$. By applying inequality \eqref{A.2} with $g^1=g$, $\xi_1=\xi$, $\zeta_1=\zeta$, $g^2=0$ , $\xi^2=0$ and  $\zeta^2=0$, we get:
\begin{equation}\label{A.5}
 (X_t^{\tau, \sigma})^2 \leq e^{\beta (T-t)} E [I(\tau, \sigma)^2|\mathcal{F}_t]+ \eta E [\int_t^T e^{\beta (s-t)}(g(s,0,0,0))^2|\mathcal{F}_t].
\end{equation}
By using the same procedure as in the proof of Proposition $\ref{oubli}$, the result follows.
\fproof

We now study the links between generalized Dynkin games (or equivalently DRBSDEs) and obstacle problems, which complete the results of this paper.

\section{Relation with partial integro-differential variational inequalities (PIDVI)}\label{sec6}
We now restrict ourselves to the Markovian case.
 Let  $b:\R \rightarrow \R$ , $\sigma :\R \rightarrow \R$ be continuous mappings, globally Lipschitz  and $\beta : \R \times E \rightarrow \R$  a measurable function  such that for some non negative real $C$, and for all $e \in E$
$$
|\beta(x,e)| \leq C \varphi (e), \quad \;\; |\beta(x,e)- \beta(x',e)| \leq C|x-x'|\varphi (e), \;\;  x, x'\in \R,$$
where $\varphi$ $\in$ $L^2_{\nu}$.
For each $(t,x) \in [0,T] \times \R$, let $( X_s^{t,x}, t \leq s \leq T)$ be the unique $\R$-valued solution of the SDE with jumps:
\begin{equation*}
X_s^{t,x}= x+ \int_t^s b(X_r^{t,x})dr+\int_t^s \sigma(X_r^{t,x})dW_r+\int_t^s\int_{{\bf E}}  \beta(X_{r^{-}}^{t,x},e) \Tilde{N}(dr,de),
\end{equation*}
and set $X_s^{t,x}=x$ for $s \leq t$.
We consider the DRBSDE associated with obstacles $\xi^{t,x}$, $\zeta^{t,x}$ of the following form:
$
\xi^{t,x}_s :=h_1(s,X_s^{t,x}),$
$\zeta^{t,x}_s :=h_2(s,X_s^{t,x}),$ $s <T$,  $\xi^{t,x}_T = \zeta^{t,x}_T:=g(X^{t,x}_T).$ We suppose that 
$g \in \mathcal{C}(\R)$, $
h_1, h_2:[0,T] \times \R \rightarrow \R
$
are jointly continuous in $t$ and $x$, and that $g$, $h_1$, $h_2$ have at most polynomial growth with respect to $x$.\\
Moreover, the obstacles $
\xi^{t,x}_s $ and
$\zeta^{t,x}_s$ are supposed to satisfy Mokobodski's condition, which holds if for example $h_1$ and $h_2$ are 
${\cal C}^{1,2}$.

We consider two functions
$\gamma$ and $f$  satisfying Assumption 2.1 in \cite{2}. More precisely,
%
we are given a map
$\gamma: \R \times E \rightarrow \R$ which is $ \mathcal{B}(\R) \otimes \mathcal{B}(E)$-measurable,  such that \\
$|  \gamma(x,e)-\gamma(x',e)| < C|x-x'|\varphi(e)$  and 
  $-1 \leq \gamma(x,e) \leq C\varphi(e)$ for each $x, x' \in \R, e \in E$.
Let $f: [0,T] \times  \R^3 \times  L_{\nu}^2  \rightarrow \R$ be a map 
supposed to be continuous in $t$ uniformly with respect to $x,y,z,k$,  and continuous in $x$ uniformly  with respect to $y,z,k$. It is also supposed to be uniformly Lipschitz with respect to $y,z,k$, and such that $f(t,x,0,0,0)$ at most polynomial growth with respect to $x$. It also satisfies that for each $t,x,y,z,k_1,k_2$,\\
$ f(t,x,y,z,k_1)- f (t,x,y,z,k_2) \geq \,  <\gamma(x,\cdot),k_1-k_2>_{\nu}$.

The driver is defined by $f(s,X_s^{t,x}(\omega),y,z,k)$.
By Theorem \ref{exiuni}, for each $(t,x) \in [0,T] \times \R$, there exists an unique solution 
 $(Y^{t,x}, Z^{t,x}, K^{t,x}, A^{t,x}, A^{' \, t,x}  ) $ of the associated DRBSDE.
Moreover, by definition, $\xi^{t,x}$ and $-\zeta^{t,x}$  are l.u.s.c. along stopping times. 
It follows that  the processes $A^{t,x}$, 
$A^{' \, t,x}$ are continuous. We define:
\begin{align}\label{4.5}
u(t,x) :=Y_t^{t,x},\quad  t \in [0,T], \;  x \in \R.
\end{align}
which is a deterministic quantity. In the following, the map $u$ is called the {\em value function} of the generalized Dynkin game. 

By the a priori estimates (see Propositions \ref{oubli} and \ref{A.4}) and the same arguments as those  used in the proofs of Lemma 3.1 and Lemma 3.2 
in \cite{2}, we derive that 
the value function $u$ is continuous in $(t,x)$ and has at most polynomial growth at infinity.

%
A solution of the obstacle problem is a function $u:[0,T] \times \R \rightarrow \R$ which satisfies the equality $u(T,x)=g(x)$ and 
%
\begin {equation}\label{4.8}
\begin{cases}
  h_1(t,x) \leq u(t,x) \leq h_2(t,x)\\
 \text{ if }u(t,x)<h_2(t,x)  \text{ then } {\cal H}u \geq 0\\
\text{  if } h_1(t,x)<u(t,x)  \text{ then } {\cal H}u \leq 0\\
\end{cases}
\end{equation}
where $L :=A+K$ 	and
\begin{itemize}
\item[$\bullet$]
$A\phi(x) := \dfrac{1}{2}\sigma^2(x)\dfrac{\p^2 \phi}{\p x^2}(x)+ b(x) \dfrac{\p \phi }{\p x}(x), \,\, B \phi(t,x)(\cdot) :=\phi(t,x+\beta(x,\cdot))-\phi(t,x) $, 
\item[$\bullet$]
$K \phi(x) :=\int_{{\bf E}} \left(\phi(x+ \beta(x,e))-\phi(x)- \dfrac{\p \phi}{\p x}(x)\beta(x,e)\right) \nu (de)$, 
\item[$\bullet$]
${\cal H}\phi(t,x):=-\dfrac{\partial \phi}{\partial t}(t,x)-L\phi(t,x)-f(t,x,\phi(t,x), (\sigma \dfrac{\partial \phi}{\p x})(t,x),B\phi(t,x))$.
\end{itemize}
Recall the classical definition of viscosity solutions.

\begin{definition}\rm
$\bullet$ A continuous function $u$ is said to be a {\em viscosity subsolution} of (\ref{4.8}) if $u(T,x)\leq g(x), x\in \R$, and if for any point $(t_0,x_0) \in [0,T) \times \R$,  we have $h_1(t_0,x_0) \leq u(t_0,x_0) \leq h_2(t_0,x_0)$ and, for any $\phi \in C^{1,2}([0,T] \times \R)$ such that $\phi(t_0,x_0)=u(t_0,x_0)$ and $\phi-u$ attains its minimum at $(t_0,x_0)$,
  if $u(t_0,x_0)>h_1(t_0,x_0)$, then $({\cal H} \phi)(t_0, x_0) \leq 0$.

$\bullet$
 A continuous function $u$  is said to be a {\em  viscosity supersolution } of (\ref{4.8}) if $u(T,x)\geq g(x), x\in \R$, and if for any point $(t_0,x_0) \in [0,T) \times \R$, we have $h_1(t_0,x_0) \leq u(t_0,x_0) \leq h_2(t_0,x_0)$ and, for any $\phi \in C^{1,2}([0,T] \times \R)$ such that $\phi(t_0,x_0)=u(t_0,x_0)$ and $\phi-u$ attains its maximum at $(t_0,x_0)$, if $u(t_0,x_0) < h_2(t_0,x_0)$ then $({\cal H}\phi)(t_0, x_0) \geq 0$.

\end{definition}
Following the same arguments as in the proof of Theorem 3.4 in \cite{2}, one can show that $u$ is viscosity subsolution of (\ref{4.8}).
By symmetry, we derive that $u$  is also a viscosity supersolution of (\ref{4.8}), which yields the following result:
\begin{theorem}
The value function $u$ is a viscosity solution (i.e. both a viscosity sub- and supersolution)  of the  obstacle problem~$(\ref{4.8})$.

\end{theorem}


In the sequel, we suppose that $E = \R^*$ and that the function $\varphi$ is defined by $\varphi(e):=1 \wedge |e|$ and is supposed to belong in $L^2_{\nu}$. We also suppose that $g$, $h_1$ and $h_2$ are bounded, and that Assumption 4.1 in \cite{2} holds. 
More precisely, \\
$ (i)   \quad  f(s,X_s^{t,x}(\omega),y,z,k) :=\overline{f}\left(s,X_s^{t,x}(\omega),y,z,\int_{\R^*}k(e)\gamma(X_s^{t,x}(\omega),e)\nu(de)\right)\textbf{1}_{s \geq t},$\\
where
$\overline{f}: [0,T] \times  \R^4  \rightarrow \R$
is a map which is continuous with respect to $t$ uniformly in $x,y,z,k$,  and continuous with respect $x$ uniformly in $y,z,k$. It is also uniformly Lipschitz with respect to $y,z,k$ and 
the map $\overline{f}(t,x,0,0,0) $ is uniformly bounded.\\
The map $k \mapsto \overline{f}(t,x,y,z,k)$ is also non-decreasing, for all $ t \in  [0,T]$, $x,y,z \in \R$.\\
$(ii) $ For each $R>0$, there exists a continuous function $m_R: \R_{+} \rightarrow \R_+$  with
$m_R(0)=0$ and \\ $|\overline{f}(t,x,v,p,q) - \overline{f}(t,y,v,p,q)| \leq m_{R}(|x-y|(1+|p|)),$   for all $ t \in  [0,T]$, $|x|, | y|\leq R, |v|\leq R,  \;p,q \in \R.$ \\
$(iii) \quad $
$|\g(x,e) - \g(y,e)| \leq C|x-y|(1\wedge e^{2})$; $0 \leq \gamma(x,e) \leq  C(1 \wedge |e|)$, $x,y \in \R, e \in \R^*.$ \\
$(iv) \quad $
 $
\overline{f}(t,x,v,p,l) - \overline{f}(t,x,u,p,l)\geq  r(u-v)$,
$u\geq v,$ $t\in [0,T]$, $x,u,v,p,l\in \R$,
where $r>0$.\\
To simplify notation, in the sequel, $\overline{f}$ is denoted by $f$.\\
The operator $B$ has now the following form: 
 $B \phi(x) :=\int_{\R^*}(\phi(x+\beta(x,e))-\phi(x)) \gamma (x,e) \nu(de)$.

%
%
%
%
%
%

\begin{theorem}[Comparison principle]\label{8.9} 
If U is a bounded viscosity subsolution  and V is a bounded viscosity supersolution of the obstacle problem (\ref{4.8}), then $U(t,x) \leq V(t,x)$,
for each $(t,x) \in [0,T] \times \R$.
\end{theorem}

\dproof For completeness, we give a sketch of proof, where we draw attention to some points which differ from the proof given in \cite{2} ( in the case of reflected BSDEs). Set
\begin{equation*}
\psi^{\epsilon, \eta}(t,s,x,y):=U(t,x)-V(s,y)- \dfrac{|x-y|^2}{\epsilon^2}-\dfrac{|t-s|^2}{\epsilon^2}- \eta^2(|x|^2+|y|^2).
\end{equation*} where $\epsilon, \eta$ are small parameters devoted to tend to $0$.
Let $M^{\epsilon, \eta}$ be a maximum of $\psi^{\epsilon, \eta}( t,s,x,y)$. This maximum is reached at some point $(t^{\epsilon, \eta},s^{\epsilon, \eta},x^{\epsilon, \eta},y^{\epsilon, \eta})$.
We define:
\begin{equation*}
\Psi_1(t,x):=V(s^{\epsilon,\eta}, y^{\epsilon, \eta})+ \dfrac{|x-y^{\epsilon, \eta}|^2}{\epsilon^2}+\dfrac{|t-s^{\epsilon, \eta}|^2}{\epsilon^2}+ \eta^2(|x|^2+|y^{\epsilon, \eta}|^2);
\end{equation*}
\begin{equation*}
\Psi_2(s,y):=U(t^{\epsilon,\eta}, x^{\epsilon, \eta})- \dfrac{|x^{\epsilon, \eta}-y|^2}{\epsilon^2}-\dfrac{|t^{\epsilon, \eta}-s|^2}{\epsilon^2}- \eta^2(|x^{\epsilon, \eta}|^2+|y|^2).
\end{equation*}

As $(t,x) \rightarrow  (U-\Psi_1)(t,x)$ reaches its maximum at $(t^{\epsilon, \eta},x^{\epsilon, \eta})$ and $U$ is a subsolution, we have the two following cases:

$\bullet$ $t^{\epsilon,\eta}=T$ and then $U(t^{\epsilon, \eta}, x^{\epsilon, \eta}) \leq f(x^{\epsilon, \eta})$, \\

$\bullet$ $t^{\epsilon,\eta} \not =T$, $h_1(t^{\epsilon, \eta}, x^{\epsilon, \eta}) \leq U(t^{\epsilon, \eta}, x^{\epsilon, \eta}) \leq h_2(t^{\epsilon, \eta}, x^{\epsilon, \eta})$ and, if $U(t^{\epsilon, \eta}, x^{\epsilon, \eta})>h_1(t^{\epsilon, \eta}, x^{\epsilon, \eta})$, we then have:
\begin{align}\label{4.16}
&-\dfrac{\p \Psi_1}{\p t}(t^{\epsilon, \eta},x^{\epsilon, \eta})-L \Psi_1(t^{\epsilon, \eta},x^{\epsilon, \eta}) - f\left(t^{\epsilon, \eta},x^{\epsilon, \eta},U(t^{\epsilon, \eta},x^{\epsilon, \eta}), (\sigma \dfrac{\p \Psi_{1}}{\p x} )(t^{\epsilon, \eta},x^{\epsilon, \eta}), B \Psi_1(t^{\epsilon, \eta},x^{\epsilon, \eta})\right) \leq 0.
\end{align}
As $(s,y) \rightarrow  (\Psi_2-V)(s,y)$ reaches its maximum at $(s^{\epsilon, \eta},y^{\epsilon, \eta})$ and $V$ is a supersolution, we have the two following cases:
\begin{itemize}
\item[$\bullet$]
 $s^{\epsilon,\eta}=T$ and $V(s^{\epsilon, \eta}, y^{\epsilon, \eta}) \geq f(y^{\epsilon, \eta})$,
\item[$\bullet$]
 $s^{\epsilon,\eta} \not =T$, $h_1(s^{\epsilon, \eta}, y^{\epsilon, \eta}) \leq V(s^{\epsilon, \eta}, y^{\epsilon, \eta}) \leq h_2(s^{\epsilon, \eta}, y^{\epsilon, \eta})$ and, if $V(s^{\epsilon, \eta}, y^{\epsilon, \eta})<h_2(s^{\epsilon, \eta}, y^{\epsilon, \eta})$ then
\begin{align*}
& -\dfrac{\p \Psi_2}{\p t}(s^{\epsilon, \eta},y^{\epsilon, \eta})-L \Psi_2(s^{\epsilon, \eta},y^{\epsilon, \eta})-f(s^{\epsilon, \eta},y^{\epsilon, \eta},V(s^{\epsilon, \eta},y^{\epsilon, \eta}),(\sigma \dfrac{\p \Psi_{2} }{\p x} )(s^{\epsilon, \eta},y^{\epsilon, \eta})), B \Psi_2(s^{\epsilon, \eta},y^{\epsilon, \eta}) \geq 0.\nonumber
\end{align*}
%
\end{itemize}

As in $\cite{2}$, we have:
$
|x^{\epsilon, \eta}-y^{\epsilon, \eta}|+|t^{\epsilon, \eta}-s^{\epsilon, \eta}| \leq C\epsilon,$ $ |x^{\epsilon, \eta}|\leq \dfrac{C}{\eta}$ and $|y^{\epsilon, \eta}| \leq \dfrac{C}{\eta}.$

Extracting a subsequence if necessary, we may suppose that for each $\eta$ the sequences $(t^{\epsilon, \eta})_{\epsilon}$ and  $(s^{\epsilon, \eta})_{\epsilon}$ converge to a common limit $t^{\eta}$, and the sequences $(x^{\epsilon, \eta})_{\epsilon}$ and $(y^{\epsilon, \eta})_{\epsilon}$ converge to a common limit $x^{\eta}.$

Here, we have to consider four cases.\\  
\textbf{1st case:} there exists a subsequence of $(t^{\eta})$ such that $t^{\eta}=T$ for all $\eta$ ( of this subsequence)\\
\textbf{2nd case:} there exists a subsequence of $(t^{\eta})$ such that $t^{\eta} \neq T$  and for all $\eta$ belonging to this subsequence, there exist a subsequence of $(x^{\epsilon,\eta})_{\epsilon}$ and a subsequence of $(t^{\epsilon,\eta})_{\epsilon}$, such that $U(t^{\epsilon, \eta},x^{\epsilon, \eta})-h_1(t^{\epsilon, \eta},x^{\epsilon, \eta}) = 0$.\\
\textbf{3rd case:} there exists a subsequence such that $t^{\eta} \not= T$, and for all $\eta$ belonging to this subsequence, there exist a subsequence of $(y^{\epsilon,\eta})_{\epsilon}$ and a subsequence of $(s^{\epsilon,\eta})_{\epsilon}$, such that $V(s^{\epsilon, \eta},y^{\epsilon, \eta})-h_2(s^{\epsilon, \eta},y^{\epsilon, \eta}) = 0.$\\
\textbf{Last case:} we are left with the case when, for a subsequence of $\eta$ we have  $t^{\eta} \not= T$, and for all $\eta$ belonging to this subsequence, there exist a subsequence of $(x^{\epsilon,\eta})_{\epsilon}$, $(y^{\epsilon,\eta})_{\epsilon}$,  $(t^{\epsilon,\eta})_{\epsilon}$ and $(s^{\epsilon,\eta})_{\epsilon}$ such that
$$
U(t^{\epsilon, \eta},x^{\epsilon, \eta})-h_1(t^{\epsilon, \eta},x^{\epsilon, \eta})>0 ; \quad
h_2(s^{\epsilon, \eta},y^{\epsilon, \eta})-V(s^{\epsilon, \eta},y^{\epsilon, \eta})>0.
$$
The first, second and fourth case are identical to  the three cases considered for reflected BSDEs (see $\cite{2}$). The third one,  which didn't appear in the case of reflected BSDEs, can be treated  similarly  to the second one.
\fproof

We derive that there exists an unique solution of the obstacle problem (\ref{4.8}) in the class of bounded continuous  functions.

\section{Appendix}

\textbf{Proof of Lemma \ref{seq}}: For completeness, we give a sketch of the proof, where we draw attention to the importance of the property $\Tilde{\xi}_T^{g} = \Tilde{\zeta}_T^{g} = 0$ a.s.\,
Set $J_{\cdot}^{(0)}=0$ and $J_{\cdot}^{\prime(0)}=0$ and  define recursively for each $n \in \mathbb{N}$, the   supermartingales:
\begin{equation}\label{sys1}
J^{(n+1)}:={\cal R}(J'^{(n)}+\Tilde{\xi}^g) \quad ; \quad J'^{(n+1)}:={\cal R} (J^{(n)}-\Tilde{\zeta}^g)
\end{equation}
 which belong to $ {\cal S}^2$. For sake of simplicity, in the above definition we have omitted the exposant $g$ in the definition of $J^{(n)}$.
Since  $\Tilde{\xi}_T^{g} = \Tilde{\zeta}_T^{g} = 0$ a.s.\,, it follows that, for each $n$,
$J^{(n)}_T = J'^{(n)}_T = 0$ a.s.
We have $J^{(0)}=0$ and $J^{'(0)}=0$. Let us prove recursively that for each $n$, $J^{'(n)}, J^{(n)}$ are well defined and nonnegative.
Suppose that $J^{'(n)}, J^{(n)}$ are well defined and nonnegative. Then $J^{(n+1)}$, $J^{'(n+1)}$ are well defined  since $(J^{'(n)}+\xi)^{-}$ and $(J^{(n)}-\zeta)^-$ belong to $\mathcal{S}^2$. Also, $J_t^{(n+1)} \geq E[J_T^{'(n)}+\Tilde{\xi}_T^g| \mathcal{F}_t]\geq 0$ a.s. since $\Tilde{\xi}_T^g=0$ a.s. Similarly, because $\Tilde{\zeta}_T^g=0$ a.s.,  $J_t^{'(n+1)} \geq 0$ a.s. By classical results, $J^{(n)}$ and $J^{'(n)}$ are RCLL supermartingales.\\
Let us prove that $J^{(n)}$ and $J^{'(n)}$ are non decreasing sequences. We have $J^{(1)} \geq 0=J^{(0)}$ and $J^{'(1)} \geq 0=J^{'(0)}$. Suppose that $J^{(n)} \geq J^{(n-1)}$ and $J^{'(n)} \geq J^{'(n-1)}$. We then have:
\begin{align}
{\cal R}(J'^{(n)}+\Tilde{\xi}^g) \geq {\cal R}(J^{'(n-1)}+\Tilde{\xi}^g) \quad ; \quad 
{\cal R}(J^{(n)}-\Tilde{\zeta}^g) \geq {\cal R}(J^{(n-1)}-\Tilde{\zeta}^g),
\end{align}
which leads to $J^{(n+1)} \geq J^{(n)}$ and $J^{'(n+1)} \geq J^{(n)}$.

Let $J^g :=\lim \uparrow J^{(n)} \; \text{ and } \; J'^g :=\lim \uparrow J'^{(n)}.$ 
Since  for each $n$,
$J^{(n)}_T = J'^{(n)}_T = 0$ a.s.\, we have $J^g_T = J^{'g}_T = 0$ a.s.\,
By classical results, $J^g$ and $J'^g$ are indistinguishable from RCLL supermartingales valued in $[0, + \infty]$, as the non decreasing limits of non negative RCLL supermartingales.
For each $n \in \mathbb{N}$, we have $J^{(n+1)}= {\cal R}(J^{'(n)}+\Tilde{\xi}^g) \leq {\cal R}(J^{'g}+\Tilde{\xi}^g)$.
Letting $n$ tend to $+\infty$, we get 
\begin{equation}\label{Snell}
J^g \leq {\cal R}(J^{'g}+\Tilde{\xi}^g).
\end{equation}
Now, for each $n \in \mathbb{N}$, $J^{(n+1)} \geq J^{'(n)}+\Tilde{\xi}^g$. By letting $n$ tend to $+\infty$, we derive that $J^g \geq J^{'g}+\Tilde\xi^g$. The characterization of ${\cal R}(J^{'g}+\Tilde\xi^g)$ yields that $J^g \geq {\cal R}(J^{'g}+\Tilde\xi^g)$. This with $\eqref{Snell}$ implies that $J^g={\cal R}(J^{'g}+\Tilde\xi^g).$ Similarly, $J^{'g}={\cal R(J^g-\Tilde \zeta^g)}$, which corresponds to the desired result \eqref{sys2}.  
We thus have $J^g \geq J^{'g} + \xi$ and 
$J^{'g} \geq J^g - \zeta$. One can prove that $J$ and $J'$ are the smallest nonnegative supermartingales 
satisfying these two inequalities (see the proof of Proposition 5.1 in \cite{KQC} for details).
\fproof
\begin{remark}\label{8}
 
We point out that the property $\tilde \xi_T^g=\tilde \zeta_T^g=0$ a.s. ensures that for each $n$,
$J^{(n)}_T = J'^{(n)}_T = 0$ a.s. We underline that if we had not made the change of variable \eqref{defz}, then $\tilde \xi^g, \tilde \zeta^g$ would be replaced  by $\ \xi, \zeta$ in the definitions of $J^{(n)}$ and  
$J^{'(n)}$. In that case, $ \xi_T= \zeta_T$ a.s. but would not necessarily be equal to $0$, and we would have $J^{(n)}_T = -J'^{(n)}_T =0$ a.s. if $n$ is even, and $\xi_T$ otherwise.
Then, the sequences $(J^{(n)}_T)_{n \in \N}$ and $(J'^{(n)}_T)_{n \in \N}$ do not converge a.s. if $P(\xi_T \neq 0). >0$ 
Also, the non negativity  property of the sequences $(J^{(n)}, n \in \N)$ and  
$(J^{'(n)}, n \in \N)$ and their non decreasing property would not necessarily hold.

\end{remark}

\textbf{Proof of Theorem \ref{f}}:\\
 We have already proved the existence. Let $(Y, Z, k, A,A')$ be a solution of the DRBSDE associated with driver process $g(t)$ and obstacles $(\xi, \zeta)$.
Let us prove that it is unique. 
%
We first show the uniqueness of $Y$. For each $S$ $\in$ $\T_0$ and for each $\varepsilon >0$, let

\begin{equation}\label{epsi}
\tau^{\varepsilon}_S:= \inf \{ t \geq S,\,\,Y_t \leq \xi_t  + \varepsilon\} \quad \sigma^{\varepsilon}_S:= \inf \{ t \geq S,\,\,Y_t \geq \zeta_t  - \varepsilon\}.
\end{equation}
Note that $\sigma^{\varepsilon}_S$  and  $\tau^{\varepsilon}_S$ $\in$ $\T_S$. Fix $\varepsilon >0$. We have that almost surely, if $t \in [S,
\tau^{\varepsilon}_S[$, then $ Y_t> \xi_t + \varepsilon$ and hence $Y_t> \xi_t$. It follows that the function $t \mapsto A^c_t$ is constant a.s. on $[S,
\tau^{\varepsilon}_S]$ and $t \mapsto A^d_t $ is constant a.s. on $[S, \tau^{\varepsilon}_S[$.
Also, $ Y_{ (\tau^{\varepsilon}_S)^-  } \geq \xi_{ (\tau^{\varepsilon}_S)^-  }+ \varepsilon\,$ a.s.\,
Since $\varepsilon >0$, it follows that $ Y_{ (\tau^{\varepsilon}_S)^-  } > \xi_{ (\tau^{\varepsilon}_S)^-  }$  a.s.\,\,, which implies that $\Delta A^d _ {\tau^{\varepsilon}_S  } =0$ a.s.\, Hence, the process $A$ is constant on $[S, \tau^{\varepsilon}_S]$.
Furthermore, by the right-continuity of $(\xi_t)$ and $(Y_t)$, we clearly have
$Y_{\tau^{\varepsilon}_S} \leq \xi_{\tau^{\varepsilon}_S} + \varepsilon \quad \mbox{a.s.}$
Similarly, one can show that the process $A'$ is constant on $[S, \sigma^{\varepsilon}_S]$ and that
$Y_{\sigma^{\varepsilon}_S} \geq \zeta_{\sigma^{\varepsilon}_S} - \varepsilon \quad \mbox{a.s.}$

Let $\tau \in \T_S $. Since $A'$ is constant on $[S, \sigma_S^{\varepsilon}]$, the process $(Y_t+\int_0^tg(s)ds,  \;  S \leq t \leq \tau \wedge \sigma_S^{\varepsilon})$ is a supermartingale. Hence
$$Y_S \geq E[Y_{\tau \wedge \sigma_S^{\varepsilon}}+ \int_S^{\tau \wedge \sigma_S^{\varepsilon}}g(s)ds  \mid \Fc_S] \quad \mbox{a.s.}$$
We also have that 
$ Y_{\tau \wedge {\sigma_S^{\varepsilon}}}= Y_{\tau}\textbf{1}_{\tau \leq \sigma_S^{\varepsilon}}+Y_{\sigma_S^{\varepsilon}}\textbf{1}_{\sigma_S^{\varepsilon}< \tau} \geq \xi_{\tau}\textbf{1}_{\tau \leq \sigma_S^*}+(\zeta_{\sigma_S^{\varepsilon}}-\varepsilon)\textbf{1}_{\sigma_S^{\varepsilon} < \tau} \text{ a.s.}$
We get\\
$
Y_S \geq E[I_S(\tau,\sigma_S^{\varepsilon})  \mid \Fc_S]-\varepsilon $ a.s.\, Similarly, one can show that for each $\sigma \in \T_S$,\\
 $Y_S \leq E[I_S(\tau_S^{\varepsilon},\sigma)  \mid \Fc_S]+\varepsilon$ a.s.\, It follows that for each $\varepsilon >0$,
$$
 ess\sup_{\tau \in \T_s}E[I_S(\tau,\sigma_S^{\varepsilon})  \mid \Fc_S]-\varepsilon \,\, \leq \,\, Y_S \,\, \leq \,\, ess\inf_{\sigma \in \T_S} E[I_S(\tau_S^{\varepsilon},\sigma)  \mid \Fc_S]+\varepsilon \,\, \text{ a.s.},
$$
that is $\overline{V}(S)-\varepsilon \,\, \leq \,\, Y_S \,\,\leq \,\, \underline{V}(S)+\varepsilon$ a.s.\, Since 
$\underline{V}(S) \leq \overline{V}(S)$  a.s. we get 
$\underline{V}(S)= Y_S = \overline{V}(S)\quad  \text{ a.s.}$
This equality holds of each stopping time $S \in \T_0$, which implies the uniqueness of $Y$. It remains to show the uniqueness of $(Z,k,A,A')$.
By the uniqueness of the decomposition of the semimartingale $Y_t+ \int_0^t g(s)ds$, there exists an unique square integrable martingale $M$ and an unique square integrable finite variation RCLL adapted process $\alpha$ with $\alpha_0=0$ such that $dY_t +g(t)dt = dM_t -d\alpha_t$. The  
martingale representation theorem applied to $M$ ensures the uniqueness of the pair $(Z,k)$ $\in$ $ \H^{2} \times  \H_\nu^{2}$.\\
The uniqueness of the processes $A$, $A'$ follows from the uniqueness of the canonical decomposition of an RCLL process with integrable variation (see Proposition \ref{canonique}).

Suppose that $A$ and $A'$ are continuous.
Since $Y$ and $\xi$ are right-continuous, we have $Y_ { \sigma^*_S } = \xi_{  \sigma^*_S}$  and $Y_ { \tau^*_S } = \zeta_{  \tau^*_S}$ a.s.
By definition of $\tau^*_S$, on $[S , \tau^*_S[$, we have $Y_t> \xi_t$ a.s. Since 
$(Y,Z,k(.),A,A')$ is the solution of the DRBSDE, $A$ is constant on $[S, \tau^*_S[$ a.s. and even on $[S, \tau^*_S]$ because $A$ is continuous. 
Similarly,  $A'$ is constant on $[S, \sigma_S^{*}]$ a.s. The process
 $(Y_t+\int_0^tg(s)ds, S \leq t \leq \tau_S^* \wedge \sigma_S^*) $ is thus a martingale. Hence, we have
$Y_S=  E[I_S(\tau_S^* ,\sigma_S^*)  \mid \Fc_S]$ a.s.\, By similar arguments as above, one can show that for each $\tau, \sigma \in \T_S$, 
$E[I_S(\tau, \sigma_S^*)  \mid \Fc_S]\leq Y_S$ and  $Y_S \leq E[I _S(\tau_S^*, \sigma)  \mid \Fc_S] $ a.s.\,, which yields that 
$(\tau_S^*, \sigma_S^*)$ is an $S$-saddle point.
\fproof

\textbf{Proof of Theorem \ref{exiuni}}:

For $\beta >0$,   $\phi \in \H^{2} $, and $l \in \H_\nu^{2}$, we introduce the norms $\| \phi \|_{\beta}^2 := E[\int_0^T e^{\beta s} \phi_s^2 ds], $ and 
 $\| l \|_{\nu,\beta}^2 := E[\int_0^T e^{\beta s} \|l_s\|_\nu^2 \, ds] $.

Let 
 $\H_{\beta,\nu}^2 $  (below simply denoted by $\H_{\beta}^2$) the space $\H^2 \times \H^2 \times \H^2_\nu$ equipped with the norm
$\| Y, Z, k(\cdot) \|_\beta^2 := \| Y \|_{\beta}^2 +  \| Z \|_{\beta}^2  +  \| k \|_{\nu,\beta}^2 $.

We define a mapping $\Phi$ from $\H_\beta^2$ into
itself as follows. Given $(U, V,l) \in \H_\beta^2$,  by Theorem \ref{f} there exists a unique process  $(Y, Z,k) = \Phi (U, V, l)$  solution of the DRBSDE  associated with driver process $g(s) = g(s, U_s , V_s, l_s)$.
Note that $(Y, Z,k)$ $\in$ $\H_\beta^2$.
Let $ A,A'$ be the associated non decreasing  processes.
Let us show that $\Phi$ is a contraction and hence admits a unique fixed point $(Y,Z,k)$ in $\H_\beta^2$, which corresponds to the unique solution of DRBSDE~\eqref{DRBSDE}.  The associated finite variation  process  is then uniquely determined in terms of $(Y,Z,k)$ and the pair $(A,A')$ corresponds to the unique canonical decomposition of this finite variation process. Let $(U^{2}, V^{2}, l^{2} )$ be another element of $\H_\beta^2$  and define $(Y^2, Z^2, k^2) = \Phi (U^2, V^2, l^2)$. Let $ A^2 , A'^2$ be the associated non decreasing processes. Set
$\overline{U} = U - U^2 , \, \,\,\overline{V} = V - V^2 , \,\,\,\overline{l} = l - l^2 $ and,
$ \overline{Y} = Y - Y^2,  \,\,\, \overline{Z} = Z - Z^2,  \,\,\, \overline{k} = k - k^2$.
By It\^o's formula, for any $\beta  > 0$, we have

\begin{align} \label{rhss}
 &\overline{Y}^2_0  + E \int^T_0 e^{\beta s } [\beta  \overline{Y}^2_s + \overline{Z}_{s}^2
 +  \| \overline k_{s}^2 \|]ds +E \sum_{0 <s \leq T} e^{\beta s }(\Delta A_s-\Delta A_s^2-\Delta A_s^{'}+\Delta A_s^{'2})^{2} \nonumber \\  
&= 2E \int^T_0 e^{\beta s}  \overline{Y}_s  [ g(s, U_s, V_s, l_s ) - g (s, U^2_s, V^2_s, l^2_s )] \, ds + 2 E[ \int^T_0 e^{\beta s} \overline{Y}_{s^-} \, dA_s - \int^T_0 e^{\beta s} \overline{Y}_{s^-}\,  dA^2_s] \nonumber \\
&- 2 E[ \int^T_0 e^{\beta s} \overline{Y}_{s^-} \, dA'_s - \int^T_0 e^{\beta s} \overline{Y}_{s^-}\,  dA'^2_s ].
\end{align}
Now, we have a.s.
$$\overline{Y}_s dA_s^c = (Y_s - \xi_s)dA_s^c - (Y_s^2 - \xi_s)dA_s^c = - (Y_s^2 - \xi_s)dA_s^c \leq 0 $$
and by symmetry, $\overline{Y}_s d{A}_s^{2c}  \geq 0$ a.s.\,
Also, we have a.s.
$$\overline{Y}_{s^-} \Delta A_s^d = (Y_{s^-} - \xi_{s^-}) \Delta A_{s}^d - (Y^2_{s^-} - \xi_{s^-}) \Delta A_{s}^d = - (Y^2_{s^-} - \xi_{s^-}) \Delta A_{s}^d \leq 0 $$
and $\overline{Y}_{s^-} \Delta {A^2}_s^d  \geq 0$ a.s.
Similarly, we have a.s.
$$\overline{Y}_s dA^{'c}_{s} = (Y_s - \zeta_s)dA_s^{'c} - (Y_s^2 - \zeta_s)dA_s^{'c} = - (Y_s^2 - \zeta_s)dA_s^{'c} \geq 0 $$
and by symmetry, $\overline{Y}_s d{A'}_s^{2c}  \leq 0$ a.s.\,
Also, we have a.s.
$$\overline{Y}_{s^-} \Delta A^{'d}_{s} = (Y_{s^-} - \zeta_{s^-}) \Delta A^{'d}_{s} - (Y^2_{s^-} - \zeta_{s^-}) \Delta A^{'d}_{s} = - (Y^2_{s^-} - \zeta_{s^-}) \Delta A^{'d}_{s} \geq 0 $$
and $\overline{Y}_{s^-} \Delta {{A'}^{2}}_s^d  \leq 0$ a.s.

Consequently, the second  and the third term  of \eqref{rhss} are non positive.
By using the Lipschitz property of $g$ and the inequality
$2 C y u \leq 2 C^2 y^2 + \frac{1}{2} u^2, $ we get 
$$\beta \|\overline{Y}\|_\beta^2 + \|\overline{Z}\|_\beta^2 + \|\overline{k}\|_{\nu, \beta}^2
\leq 6 C^2 \|\overline{Y}\|_\beta^2 + \frac{1}{2} (
\|\overline{U}\|_\beta^2 + \|\overline{V}\|_\beta^2 + \|\overline{l}\|_{\nu, \beta}^2
).$$
Choosing $\beta  = 6 C^2 + 1$, we deduce
$ \| (\overline{Y}, \overline{Z}, \overline{k}) \|_\beta^2 \leq \frac{1}{2}
\| (\overline{U}, \overline{V}, \overline{l}) \|_\beta^2.  $\\
The last assertion of the theorem follows from Theorem \ref{continuous}.
\fproof

\begin{lemma} \label{diff}
Let $g$ be a Lipschitz driver. Suppose that $g$ is differentiable (or equivalently Fr\'echet-differentiable) of class ${\cal C}^1$ with respect to $l$    and such that $ \nabla_{l} g (t,\omega,y, z,l)$ satisfies that for each $(t, \omega,y, z,l)$,  $$\nabla_{l} g (t,\omega,y, z,l)(e)\geq -1\,\,\,  \;\; \text{ and }
\,\,  \;\; \vert \nabla_{l} g (t,\omega, y,z,l)(e) \vert \leq\psi (e), \quad \forall e \in E$$ where $\psi$ $\in$ $L^2_{\nu}$.\\ 
Then $g$ satisfies Assumption \ref{Royer}.
\end{lemma}

\dproof 
Let $(Y,Z,l^1,l^2)$ in $S^{2} \times \H^{2} \times (\H_{\nu}^{2})^2$.
Let $(t, \omega) \in [0,T] \times \Omega$. By Lagrange's theorem, there exists $\lambda \in [0,1]$ such that 
\begin{align*}
g(t, \omega,Y_t(\omega), Z_t(\omega), l^1_t(\omega)) &- g(t,\omega,Y_t(\omega),Z_t(\omega),l^2_t(\omega)) \\
&=  <\nabla_l g(t, \omega,Y_t(\omega), Z_t(\omega), l^2_t(\omega)+\lambda(l^1_t(\omega)-l^2_t(\omega)), l^1_t(\omega)-l^2_t(\omega)>_{\nu}
\end{align*}
By the section theorem, 
there exists a predictable process $(\lambda_t)$ valued in $[0,1]$ such that 
\begin{equation*}\label{autre}
g(t, Y_t, Z_t, l^1_t) - g(t,Y_t,Z_t,l^2_t) = \langle \gamma_t\,,\, l^1_t- l^2_t \rangle_\nu , \;\; t \in [0,T],\; \; dt\otimes dP \text{ a.s.}.
\end{equation*}
where $\gamma_t:= \nabla_l g(t, Y_t, Z_t, l^2_t+\lambda_t(l^1_t-l^2_t))$.
Hence, $g$ satisfies Assumption \ref{Royer}.
\fproof

%
%

We now easily show an ${\cal E}$- Doob-Meyer decomposition of ${\cal E}$-supermartingales , which generalizes the results given in \cite{Peng} and \cite{R}  under stronger assumptions. Moreover, our proof gives an alternative proof of these previous results.
\begin{definition}\label{defmart}
Let $Y \in \cal S^2$ . The process $(Y_t)$ is said to be a strong ${\cal E}$-supermartingale (resp ${\cal E}$-submartingale), if ${\cal E}_{\sigma ,\tau}(Y_{\tau}) \leq Y_{\sigma}$ (resp. ${\cal E}_{\sigma ,\tau}(Y_{\tau})\geq Y_{\sigma})$ a.s.\, on $\sigma \leq \tau$,  for all $ \sigma, \tau \in \T_0$. 
\end{definition}

\begin{proposition}\label{compref} Suppose that $g$ satisfies Assumption $\eqref{Royer}$. 
\begin{itemize}
\item
 Let $A$  be a non decreasing (resp non increasing) RCLL predictable process in ${\cal S}^2$ with $A_0=0$. Let $(Y,Z,k) \in \cal S^2 \times \mathbb{H}^2 \times \mathbb{H}^2_{\nu}$  following the dynamics:
 \begin{eqnarray}\label{dy}
 -dY_{s}= g(s,Y_{s}, Z_{s}, k_s)ds + dA_s - Z_{s}dW_{s} -\int_{{\bf E}}  k_s(e) \tilde{N}(ds,de).
 \end{eqnarray}
Then the process $(Y_t)$ a strong ${\cal E}$-supermartingale (resp ${\cal E}$-submartingale).
\item 
(${\cal E}$-Doob-Meyer decomposition)
Let $(Y_t)$ be a strong ${\cal E}$-supermartingale (resp. ${\cal E}$-submartingale).
Then, there exists a non decreasing (resp non increasing) RCLL predictable process $A$ in ${\cal S}^2$ with $A_0=0$ and 
  $(Z,k) \in   \H^2 \times \H^2_{\nu}$  such that (\ref{dy}) holds.
\end{itemize}
\end{proposition}

\dproof Suppose $A$ is non decreasing. Let $(X^{\tau},\pi^{\tau}, l^{\tau})$ be the solution of the BSDE associated with driver $g$, terminal time $\tau$, and terminal condition $Y_{\tau}$, 
%
Since $g$ satisfies Assumption \ref{Royer} and since  $g(s,y,z ,k)ds + dA_{s} \geq g(s,y,z,k)ds $,  the comparison theorem for BSDEs (see Theorem 4.2 in  \cite{16})  gives that $Y_\sigma \geq X^{\tau}_\sigma = {\cal E}_{\sigma , \tau}(Y_{\tau}) $ a.s. on $\{\sigma  \leq \tau\}$. The case when $A$ is non-increasing can be shown similarly.

Let us show the second assertion.  Fix $S$ $\in$ $\T_0$.
Since $(Y_t)$ is a strong ${\cal E}$-supermartingale, we derive that for each $\tau \in \T_S$, we have $Y_{S} \geq {\cal E}_{S ,\tau}(Y_{\tau})$ a.s.\, We thus get
  $$Y_S  \geq  {\rm ess} \sup_{\tau \in \T_S} {\cal E}_{S ,\tau}(Y_{\tau})  \quad a.s.$$
  Now, by definition of the essential supremum, $Y_S  \leq  {\rm ess} \sup_{\tau \in \T_S} {\cal E}_{S ,\tau}(Y_{\tau}) $ a.s.\, because $S \in \T_S$. 
 The two above inequalities imply that 
  $$ Y_S  =  {\rm ess} \sup_{\tau \in \T_S} {\cal E}_{S ,\tau}(Y_{\tau})  \quad a.s.$$
  By the characterization theorem  (Theorem 3.3 in \cite{17}) of the solution of a reflected BSDE (associated with an obstacle  supposed to be only RCLL), the process $(Y_t)$ coincides with the solution of the reflected BSDE associated with the RCLL obstacle $(Y_t)$. The result follows.
\fproof

We now show the following result on RCLL adapted processes with integrable total variation.
\begin{proposition}\label{canonique} 
Let $(\Omega, {\cal F}, P)$ be a probability space equipped with a completed\\
 right-continuous filtration $({\cal F}_t)_
{0 \leq t \leq T}$.
Let $\alpha=(\alpha_t)_{0 \leq \leq T}$ be a RCLL adapted process with integrable total variation, that is, $E(\int_0^T \vert d\alpha_t \vert) < \infty$.\\
 There exists an unique  pair $(A, A') \in (\mathcal{A}^1)^2 $ such that $\alpha=A-A'$ with $dA_t \perp dA'_t$.\\
 This decomposition is called the {\em canonical decomposition} of the process $\alpha$.\\
 Moreover, if $(B, B') \in (\mathcal{A}^1)^2 $ satisfies  $\alpha = B-B'$, then $dA_t << dB_t$ in the (probabilistic) sense, that is, for each $K \in {\cal P}$ with $\int_0^T {\bf 1}_K dB_t = 0$ a.s.\,, then $\int_0^T {\bf 1}_K dA_t = 0$ a.s.
\end{proposition}

\dproof
By classical results, the process $\alpha$ can be written as $\alpha = B-B'$ with 
$B, B'$ $\in \mathcal{A}^1 $.
Let $C_t:=B_t+B'_t$. This process belongs to $\mathcal{A}^1$. For almost every $\omega$, the measures $dB_\cdot(\omega)$ and $dB'_\cdot(\omega)$ on $[0,T]$ are absolutely continuous with respect to $dC_\cdot(\omega)$. By using the Radon-Nikodym Theorem for predictable RCLL non decreasing processes (see Th. 67, Chap. VI in \cite{DM2}),  there exist non negative predictable processes $H$ and $H'$ such that for each 
$t \in [0,T]$, $$B_t=\int_0^t H_s dC_s \quad {\rm and} \quad B'_t=\int_0^t H'_s dC_s\quad {\rm a.s}$$ 
Let $A$ and $A'$ be the processes defined by  $$A_t:=\int_0^t (H_s-H'_s)^+dC_s\quad {\rm and} \quad A'_t:=\int_0^t (H_s-H'_s)^-dC_s.$$
They belong to $\mathcal{A}^1 $. 
Now, the set $D$ defined by 
$$D:=\{ (t,\omega)\,,\,H_t(\omega)-H'_t (\omega)\geq 0 \} $$ 
belongs to $ \mathcal{P}$. 
We have $\int_0^T \textbf{1}_{D_t^c} dA_t =\int_0^T  \textbf{1}_{\{H_t-H'_t< 0 \} }
 (H_t-H'_t)^+dC_t =0$ a.s.\,
Similarly  $\int_0^T \textbf{1}_{D_t} dA'_t=0$ a.s.\,, which implies  that $dA_t \perp dA'_t.$\\ 
It remains to show the uniqueness of this  decomposition. 
Since $dA_t \perp dA'_t$, it follows that, for almost every $\omega$, the deterministic measures $dA_t (\omega)$ and $dA'_t(\omega)$
are mutually singular in the classical analysis sense.
Hence, for almost every $\omega$, the non decreasing maps $A_.(\omega)$ and $A'_.(\omega)$  correspond to the unique canonical decomposition of the RCLL bounded variational map $\alpha.(\omega)$ by a well-known analysis result. This implies  the uniqueness of $A$, $A'$.\\
Moreover, since $(H_t-H'_t)^+ \leq H_t$, the last assertion holds.
\fproof

\begin{remark}\label{equi}
Note that it is obvious that, if $dA_t$ and $dA'_t$ are mutually singular in the probabilistic sense (see Definition \ref{proba}), then for almost every $\omega$, the deterministic measures on $[0,T]$ $dA_t(\omega)$ and $dA'_t(\omega)$ are mutually singular in the classical analysis sense. The converse is not so immediate. However, it holds by the above property. 
\end{remark}

\end{document}